\documentclass[12pt,reqno]{amsart}

\title[]{Invariant Measures for a Stochastic Electroconvection Model}

\author{Elie Abdo}
\address{Department of Mathematics, Temple University, Philadelphia, PA 19122}
\email{abdo@temple.edu}
\author{Mihaela Ignatova}
\address{Department of Mathematics, Temple University, Philadelphia, PA 19122}
\email{ignatova@temple.edu}

\usepackage[margin=1in]{geometry}
\usepackage{amsmath, amsthm, amssymb}
\usepackage{times}
\usepackage{color}
\usepackage{hyperref}

\newcommand{\la}{\label}
\newcommand{\fr}{\frac}
\newcommand{\na}{\nabla}
\newcommand{\be}{\begin{equation}}
\newcommand{\ee}{\end{equation}}
\newcommand{\ba}{\begin{array}{l}}
\newcommand{\ea}{\end{array}}

\newtheorem{thm}{Theorem}

\newcommand{\beg}{\begin}

\renewcommand{\l}{\Lambda}

\usepackage{MnSymbol,wasysym}

\def\TT{{\mathbb T}}

\def\PP{\mathbb P}
\def\d{\mathrm{\textbf{d}}}
\def\E{\mathrm{\textbf{E}}}

\date{today}
\begin{document}
\begin{abstract} 
We consider a stochastic electroconvection model describing the nonlinear evolution of a surface charge density in a two-dimensional fluid with additive stochastic forcing.  We prove the existence and uniqueness of solutions and we show that the corresponding Markov semigroup is weak Feller.  We also prove the existence of invariant measures for the Markov transition kernels associated with the model.   
\end{abstract} 

\vspace{.5cm}

\maketitle
\section{Introduction} \la{intro}

We consider a stochastic electroconvection model describing the evolution of a surface charge density interacting with a two-dimensional fluid. The surface charge density $q$ evolves according to the stochastic partial differential equation 
\be  \la{intro1}
\d q  + \nabla \cdot J dt = \tilde{g} dW.
\ee The current density $J$ is given by
\be 
J = E + qu 
\ee where 
\be 
E = - \na \Phi - \na \l^{-1} q,
\ee
and $\Phi$ is a potential due to applied voltage restricted to the surface whereas $\l^{-1} q$ is the potential due to the surface charge density $q$ restricted to the surface. Here $\Lambda$ denotes the square root of the two-dimensional periodic Laplacian, and $\l^{-1}$ denotes its inverse.  
The fluid velocity $u$ obeys a stochastic forced Navier-Stokes equation given by 
\be 
\d u + u \cdot \nabla u dt - \Delta u dt + \nabla p dt = qE dt + f dt + gdW, 
\ee 
and the divergence-free condition
\be 
\nabla \cdot u = 0,
\ee where $f$ are body forces and $p$ is the fluid pressure. The potential $\Phi$ is assumed to be time independent and smooth whereas the body forces $f$ are assumed to be time independent and divergence-free. We denote by
$W(t,w) = (W_1, ..., W_n)$ a collection of standard independent Brownian motions. The stochastic noise processes $gdW$ and $\tilde{g} dW$ are given by
\be 
gdW = \sum\limits_{l=1}^{n} g_l (x) dW_l (t,w)
\ee and 
\be  \la{intro2}
\tilde{g} dW =  \sum\limits_{l=1}^{n} \tilde{g}_l (x) dW_l (t,w),
\ee  where $g = (g_1, ..., g_n)$ and $\tilde{g} = (\tilde{g}_1, ..., \tilde{g}_n)$ are time-independent and the components of $g$ are divergence-free.
The system of equations \eqref{intro1}--\eqref{intro2} is posed on the two-dimensional torus $\mathbb{T}^2 = [0, 2\pi]^2$ with periodic boundary conditions.

In \cite{AI}, we considered the two-dimensional periodic deterministic electroconvection model \eqref{intro1}--\eqref{intro2}, where the equations are not forced by noise, and established the existence and uniqueness of global regular solutions, provided that the initial data is sufficiently regular. We addressed the long-time behavior of solutions and proved the existence of a finite-dimensional global attractor. 
In \cite{ceiv}, global existence of regular solutions of the deterministic model \eqref{intro1}--\eqref{intro2} with homogeneous Dirichlet boundary conditions was established in the absence of body forces in the fluid $(f=0)$.  

In this paper, we study the stochastic model described by \eqref{intro1}--\eqref{intro2} in the presence of noises forcing the equations satisfied by the charge density $q$ and the velocity $u$. We show that the stochastic system \eqref{intro1}--\eqref{intro2} has unique global solutions when the initial deterministic charge density is at least $L^4(\TT^2)$ regular and the initial deterministic velocity is at least $H^1(\TT^2)$ regular. The existence of solutions is obtained by taking a mollification of \eqref{intro1}--\eqref{intro2}, establishing uniform bounds for the mollified solutions, and using the Banach Alaoglu theorem in order to obtain weak convergence. The identification of the drift in the case of the stochastic electroconvection model \eqref{intro1}--\eqref{intro2} is highly challenging. The reason is that the nonlinearity $qRq$ is not weakly continuous in the spaces we have control in. The remedy is a coercive estimate \eqref{contidea} and use of ideas from \cite{MS} where the authors proved the existence of global solutions for the forced stochastic Navier-Stokes equations. As a consequence of the existence result, the Markov transition kernels are defined for data $q_0 \in L^4$ and $u_0 \in H^1$. 

In the absence of the potential $\Phi$ (that is $\Phi = 0$), we prove that the stochastic model \eqref{intro1}--\eqref{intro2} has an invariant measure. The requirement of a vanishing potential $\Phi$ is due to the fact that the velocity $u$ does not maintain a zero spatial average for all positive times regardless of whether or not the average of the initial velocity vanishes. The term $q\na \Phi$ forcing the velocity equation does not have a zero mean over $\TT^2$, and so the expectation of the $L^2$ norm of the velocity might grow exponentially in time. The Krylov-Bogoliubov procedure is applied to prove the existence of an invariant measure after we obtain  bounds 
\be \la{regcr}
\frac{1}{t}\E\int_{0}^{t} (\|q(s)\|_{H^s}^2 + \|u(s)\|_{H^2}^2) ds \le C
\ee 
for some $s > 1/2$ when $u_0 = q_0 = 0$ (The Krylov-Bogoliubov procedure is then applicable because of the compactness of $H^s \oplus H^2$ in $L^4 \oplus H^1$.) The required $H^s(\TT^2)$ regularity is difficult to obtain, due to the nonlinear terms involved in \eqref{intro1}--\eqref{intro2} together with the insufficient critical regularity obtained from the dissipative term. By contrast, in the subcritical case where the term $\l q$ in the surface charge density equation is replaced by $\l^{\alpha}q$ for some $\alpha > 1$, the desired bounds  \eqref{regcr} are directly obtained due to the higher regularity of the dissipation.

The existence of ergodic invariant measures for stochastic partial differential equations has been extensively studied. The existence of an invariant measure for the stochastic Navier-Stokes equations was obtained in \cite{F}. In \cite{DX}, global existence and uniqueness of strong solutions for the 2D stochastic Navier-Stokes equations on the two-dimensional torus was established and existence of invariant measures was obtained on the base of the Krylov-Bogoliubov averaging procedure. In \cite{GKVZ}, the authors proved existence of invariant measures for the 3D stochastic primitive equations by establishing moment bounds for strong solutions. In \cite{CGV}, existence and uniqueness of an ergodic invariant measure was obtained for the 2D fractionally dissipated periodic stochastic Euler equation by deriving moment bounds in Sobolev spaces that grow linearly in time.

This paper is organized as follows. In section \ref{WeakSol}, we prove that the system \eqref{intro1}--\eqref{intro2} has a unique global solution provided that the initial charge density has a zero spatial average and is $L^4$ integrable, the initial velocity is divergence-free and is weakly differentiable, and the noise is sufficiently regular. Then we define the semigroup associated with \eqref{intro1}--\eqref{intro2} in section \ref{sec3} and we prove that it is weak Feller. In the absence of potential $(\Phi = 0)$, we prove in section \ref{sec5} the existence of an invariant measure for the Markov transition kernels associated with the electroconvection model \eqref{intro1}--\eqref{intro2} based on the Krylov-Bogoliubov averaging procedure. Finally, we treat the stochastic subcritical case in section \ref{sec8} and we obtain the existence of an invariant measure.

\section{Existence and Uniqueness of Solutions} \la{WeakSol}

Let $(\Omega, \mathcal{F}, P)$ be a probability space. Let $W(t,w) = (W_1, ..., W_n)$ be a collection of independent standard Brownian motions.   
Let $T > 0$. We consider the It\^o stochastic model 
\be \begin{cases} \la{stochastic}
\d q + u \cdot \na q dt + \l q dt = \Delta \Phi dt + \sum\limits_{l=1}^{n} \tilde{g}_l dW_l 
\\\d u + u \cdot \na u dt - \Delta u dt + \na p dt = - q Rq dt - q \na \Phi dt + fdt + \sum\limits_{l=1}^{n} g_l dW_l
\\ \na \cdot u = 0
\end{cases}
\ee on $\mathbb{T}^2 \times [0,T] \times \Omega$, with initial data $q(x,0) = q_0$ and $u(x,0) = u_0$. 
The unknowns $q(x, t, w)$, $u(x,t,w) = (u_1(x,t,w), u_2 (x,t,w))$, and $p (x,t,w)$ depend on three different variables: position $x \in \mathbb{T}^2$, time  $t \in [0,T]$, and outcome  $w \in \Omega$. The body forces $f$ and the potential $\Phi$ depend only on the position variable $x$. The forces $f$ are smooth, divergence-free and have a zero space average. The potential $\Phi$ is assumed to be smooth. The functions $\tilde{g}_l(x)$ and $g_l(x)$ are assumed to be time-independent and square-integrable over the torus $\TT^2$. The functions $g_l$ are assumed to be divergence-free and the functions $\tilde{g}_l$ are assumed to have mean zero for all $l \in \left\{1, ..., n\right\}$. Here $\l$ is the periodic fractional Laplacian of order one and $R = (R_1, R_2)$ is the periodic Riesz transform.

We show the existence of solutions for the stochastic system \eqref{stochastic}.
For each $\epsilon \in (0,1]$, we let $J_{\epsilon}$ be the standard mollifier operator, and we let $(q^{\epsilon}, u^{\epsilon})$ be the solution of the stochastic system 
\be 
\begin{cases} \la{nonlinear}
\d q^{\epsilon} + u^{\epsilon} \cdot \na q^{\epsilon}  dt + \l q^{\epsilon} dt -\epsilon \Delta q^{\epsilon}
= \Delta \Phi dt + \sum\limits_{l=1}^{n} J_{\epsilon} \tilde{g}_l dW_l
\\ \d u^{\epsilon} + u^{\epsilon} \cdot \na u^{\epsilon} dt - \Delta u^{\epsilon} dt + \na p^{\epsilon} dt
= - q^{\epsilon} Rq^{\epsilon}  dt - q^{\epsilon} \na \Phi dt + f dt 
+ \sum\limits_{l=1}^{n} J_{\epsilon} g_l dW_l
\\\na \cdot u^{\epsilon} = 0
\end{cases}
\ee
with smoothed out initial data $q_0^{\epsilon} = J_{\epsilon}q_0, u_0^{\epsilon} = J_{\epsilon}u_0.$

\beg{prop} \la{prop1} Let $\epsilon \in (0,1]$ and let $T > 0$. Let $q_0 \in L^2$ have mean zero over $\TT^2$. Let $u_0 \in L^2$ be divergence-free. The stochastic system \eqref{nonlinear} has a solution $(q^{\epsilon}, u^{\epsilon})$ on $[0,T]$ such that $q^{\epsilon}$ has mean zero, $u^{\epsilon}$ is divergence-free, and $(q^{\epsilon}, u^{\epsilon})$ satisfies the following properties:
\begin{enumerate}
\item[(i)] If $\tilde{g}_{l} \in L^2(\TT^2)$ for all $l \in \left\{1, ..., n\right\}$, then $q^{\epsilon}$ is uniformly bounded in 
\be 
L^2(\Omega; L^{\infty}(0,T; L^2(\TT^2))) \cap L^2(\Omega; L^2(0,T; H^{\fr{1}{2}}(\TT^2)))
\ee
and satisfies 
\beg{align} \la{SS1}
\E \left\{\sup\limits_{0 \le t \le T} \|q^{\epsilon}\|_{L^2}^2  + 2\int_{0}^{T} \|\l^{\fr{1}{2}} q^{\epsilon} \|_{L^2}^2 ds  \right\}  
\le 4\|q_0\|_{L^2}^2 + C\left(\|\l^{\fr{3}{2}} \Phi\|_{L^2}^2 + \|\tilde{g}\|_{L^2}^2 \right)T. 
\end{align}

\item[(ii)] Let $p \in [4, \infty)$. If $\tilde{g}_{l} \in L^2(\TT^2)$ for all $l \in \left\{1, ..., n\right\}$, then $q^{\epsilon}$ is uniformly bounded in 
\be 
L^p (\Omega ; L^{\infty}(0,T; L^2(\TT^2)))
\ee and satisfies
\beg{align} \la{SS2}
\E \left(\sup\limits_{0 \le t \le T} \|q^{\epsilon}\|_{L^2}^p\right) 
&+ \fr{p^2}{2} \E \left(\int_{0}^{T} \|q^{\epsilon}\|_{L^2}^{p-2} \|\l^{\fr{1}{2}} q^{\epsilon}\|_{L^2}^2 ds \right) \nonumber
\\&\quad\quad \le 2p\|q_0\|_{L^2}^p  
+ C_p \left(\|\Delta \Phi\|_{L^2}^p  + \|\tilde{g}\|_{L^2}^p \right) T
+ C_p \|\tilde{g}\|_{L^2}^p T^{\fr{p}{2}}.
\end{align}

\item[(iii)] If If $\tilde{g}_{l} \in L^2(\TT^2)$ and $g_l \in L^2(\TT^2)$ for all $l \in \left\{1, ..., n\right\}$, then $u^{\epsilon}$ is uniformly bounded in 
\be 
L^2(\Omega; L^{\infty}(0,T;L^2(\TT^2))) \cap L^2(\Omega; L^2(0,T; H^1(\TT^2)))
\ee and satisfies
\beg{align} \la{SS3}
&\E \left\{\sup\limits_{0 \le t \le T} \|u^{\epsilon}\|_{L^2}^2 
+ \int_{0}^{T} \|\na u^{\epsilon}\|_{L^2}^2 dt \right\} 
\le  C(\|u_0\|_{L^2}, \|q_0\|_{L^4}, f, \Phi, \tilde{g}, g)e^{4T}.
\end{align}

\item[(iv)] If $q_0 \in L^4(\TT^2)$ and $\tilde{g}_{l} \in L^4(\TT^2)$ for all $l \in \left\{1, ..., n\right\}$, then $q^{\epsilon}$ is uniformly bounded in 
\be 
L^4 (\Omega ; L^{\infty}(0,T ; L^4(\TT^2))) 
\ee and satisfies 
\beg{align}  \la{SS4}
\E \left\{ \sup\limits_{0 \le t \le T} \|q^{\epsilon}\|_{L^4}^4\right\} 
&+ 4c \E \left\{\int_{0}^{T} \|q^{\epsilon}\|_{L^4}^4 dt \right\}
\le 8\|q_0\|_{L^4}^4 \nonumber\\
&\quad\quad+ C\|\Delta \Phi\|_{L^4}^4 T +  C \left(\sum\limits_{l=1}^{n} \|\tilde{g}_l\|_{L^4}^2\right)^2 T
+ C \left(\sum\limits_{l=1}^{n} \|\tilde{g}_l\|_{L^4}^2 \right)^2T^2.
\end{align}

\item[(v)] Let $p \in [8, \infty)$ be an even integer. If $q_0 \in L^4(\TT^2)$ and $\tilde{g}_{l} \in L^4(\TT^2)$ for all $l \in \left\{1, ..., n\right\}$, then $q^{\epsilon}$ is uniformly bounded in 
\be 
L^p(\Omega ; L^{\infty}(0,T; L^4(\TT^2)))
\ee
and satisfies 
\beg{align}  \la{SS5}
\E \left\{\sup\limits_{0 \le t \le T} \|q^{\epsilon}\|_{L^4}^p \right\}
+ \fr{cp^2}{2} \left\{\int_{0}^{T} \|q^{\epsilon}\|_{L^4}^p \right\}
&\le 2p\|q_0\|_{L^4}^p + C_p \|\Delta \Phi\|_{L^4}^pT \nonumber\\
&\quad\quad+ C_p\left(\sum\limits_{l=1}^{n} \|\tilde{g}_l\|_{L^4}^2 \right)^{\fr{p}{2}} T 
+ C_p  \left(\sum\limits_{l=1}^{n} \|\tilde{g}_l\|_{L^4}^2 \right)^{\fr{p}{2}}T^{\fr{p}{2}}.
\end{align}

\item[(vi)] Let $p \ge 4$.  If $q_0 \in L^4(\TT^2)$, $\tilde{g}_{l} \in L^4(\TT^2)$, and $g_l \in L^2(\TT^2)$ for all $l \in \left\{1, ..., n\right\}$, then $u^{\epsilon}$ is uniformly bounded in  
\be 
L^{p} (\Omega; L^{\infty}(0,T ; L^2(\TT^2)))
\ee
and satisfies 
\beg{align} \la{SS6} 
&\E \left\{\sup\limits_{0 \le t \le T} \|u^{\epsilon}\|_{L^2}^p \right\} 
+ \E \left\{\int_{0}^{T} \|u^{\epsilon} \|_{L^2}^{p-2} \|\na u^{\epsilon}\|_{L^2}^2  dt \right\}  
\le C(p, \|q_0\|_{L^4}, \|u_0\|_{L^2}, f, \Phi, g, \tilde{g})e^{pT}.
\end{align}

\item[(vii)] If $q_0 \in L^4(\TT^2)$, $u_0 \in H^1(\TT^2)$, $\tilde{g}_{l} \in L^4(\TT^2)$, and $g_l \in H^1(\TT^2)$ for all $l \in \left\{1, ..., n\right\}$, then $u^{\epsilon}$ is uniformly bounded in  
\be 
L^{2} (\Omega; L^{\infty}(0,T ; H^1(\TT^2))) \cap L^2(\Omega; L^2(0, T; H^2(\TT^2)))
\ee
and satisfies
\beg{align} \la{SS7} 
\E  \left\{\sup\limits_{0 \le t \le T} \|\na u^{\epsilon}(t)\|_{L^2}^2 \right\}
&+ E\left\{\int_{0}^{T} \|\Delta u^{\epsilon} (s)\|_{L^2}^{2}  ds \right\} \nonumber
\\&\quad \le C(\|\na u_0\|_{L^2}, \|q_0\|_{L^4}) + C(\Phi, f, g, \tilde{g})T + C(\tilde{g})T^2.
\end{align} 
\end{enumerate}
\end{prop}

For simplicity, we ignore the viscous term $-\epsilon \Delta q^{\epsilon}$ in the proof of proposition~\ref{prop1} below because it does not have any major contribution in estimating the solutions of the mollified system \eqref{nonlinear} and vanishes as we take the limit $\epsilon\to 0$.

\textbf{Proof of (i).} We apply It\^o's lemma pointwise in $x$ to the stochastic process $F(X_t(w))$ where $F(\xi) = \xi^2$ and $X_t(w) = q^{\epsilon}$, and we obtain  
\be
\d (q^{\epsilon})^2 
= -2 q^{\epsilon} (u^{\epsilon} \cdot \na q^{\epsilon}) dt 
-2 q^{\epsilon} \l q^{\epsilon} dt 
+2  q^{\epsilon}  \Delta \Phi dt 
+ \sum\limits_{l=1}^{n} (J_{\epsilon} \tilde{g}_l)^2 dt
+ 2 \sum\limits_{l=1}^{n} q^{\epsilon} J_{\epsilon} \tilde{g}_l dW_l.
\ee 
Next we integrate in the space variable over $\TT^2$. 
In view of the divergence-free condition obeyed by $u^{\epsilon}$, the nonlinear term vanishes, that is 
\be 
(u^{\epsilon} \cdot \na q^{\epsilon}, q^{\epsilon})_{L^2} = 0,
\ee  which yields the energy equality
\be \la{ql2}
\d \|q^{\epsilon}\|_{L^2}^2 + 2\|\l^{\fr{1}{2}} q^{\epsilon}\|_{L^2}^2
= 2 (\Delta \Phi, q^{\epsilon})_{L^2} + \sum\limits_{l=1}^{n} \|J_{\epsilon} \tilde{g}_l\|_{L^2}^2 dt 
+ 2 \sum\limits_{l=1}^{n} (J_{\epsilon} \tilde{g}_l, q^{\epsilon})_{L^2} d W_l.
\ee We estimate 
\be 
|(\Delta \Phi, q^{\epsilon})_{L^2}| = |(\l^{\fr{3}{2}} \Phi, \l^{\fr{1}{2}} q^{\epsilon})_{L^2}| \le \fr{1}{2} \|\l^{\fr{3}{2}} \Phi\|_{L^2}^2 + \fr{1}{2} \|\l^{\fr{1}{2}} q^{\epsilon}\|_{L^2}^2
\ee using the H\"older and Young inequalities. We obtain the differential inequality 
\be 
\d \|q^{\epsilon}\|_{L^2}^2 + \|\l^{\fr{1}{2}} q^{\epsilon}\|_{L^2}^2 dt
\le \|\l^{\fr{3}{2}} \Phi\|_{L^2}^2 dt +\|\tilde{g}\|_{L^2}^2 dt + 2\sum\limits_{l=1}^{n} (J_{\epsilon} \tilde{g}_l, q^{\epsilon})_{L^2} dW_l.
\ee
Integrating in time from $0$ to $t$, we get 
\beg{align} \la{S1}
&\|q^{\epsilon}(t,w)\|_{L^2}^2 + \int_{0}^{t} \|\l^{\fr{1}{2}} q^{\epsilon} (s,w)\|_{L^2}^2 ds \nonumber
\\&\quad\le \|q_0\|_{L^2}^2 + \left(\|\l^{\fr{3}{2}} \Phi\|_{L^2}^2 + \|\tilde{g}\|_{L^2}^2  \right)t 
+ 2 \int_{0}^{t} \sum\limits_{l=1}^{n} (J_{\epsilon} \tilde{g}_l, q^{\epsilon})_{L^2} dW_l. 
\end{align}
We take the supremum over all $t \in [0,T]$,
\beg{align} 
&\sup\limits_{0 \le t \le T} \|q^{\epsilon}(w)\|_{L^2}^2 + \int_{0}^{T} \|\l^{\fr{1}{2}} q^{\epsilon} (s,w)\|_{L^2}^2 ds  \nonumber
\\&\quad\le 2\|q_0\|_{L^2}^2 + 2\left(\|\l^{\fr{3}{2}} \Phi\|_{L^2}^2 + \|\tilde{g}\|_{L^2}^2 \right)T 
+ 4 \sup\limits_{0 \le t \le T} \left|\int_{0}^{t} \sum\limits_{l=1}^{n} (J_{\epsilon} \tilde{g}_l, q^{\epsilon})_{L^2} dW_l\right|. 
\end{align} Now we apply the expectation $\E$. In view of the martingale estimate (see \cite{DPZ}), 
\be 
\E \left\{\sup\limits_{0 \le t \le T} \left|\int_{0}^{t} \sum\limits_{l=1}^{n} (J_{\epsilon} \tilde{g}_l, q^{\epsilon})_{L^2} dW_l \right| \right\}
\le C \E \left\{\left(\int_{0}^{T} \sum\limits_{l=1}^{n} (J_{\epsilon} \tilde{g}_l, q^{\epsilon})_{L^2}^2 dt \right)^{\fr{1}{2}}\right\},
\ee we have 
\beg{align}
&\E \left\{\sup\limits_{0 \le t \le T} \left|\int_{0}^{t} \sum\limits_{l=1}^{n} (J_{\epsilon} \tilde{g}_l, q^{\epsilon})_{L^2} dW_l \right| \right\}
\le C \E \left\{\left(\int_{0}^{T} \|q^{\epsilon}\|_{L^2}^2\|\tilde{g}\|_{L^2}^2 dt \right)^{\fr{1}{2}} \right\} \nonumber
\\&\quad\le \E \left\{\left(\sup\limits_{0 \le t \le T} \|q^{\epsilon}\|_{L^2} \right) \left(C\int_{0}^{T} \|\tilde{g}\|_{L^2}^2 dt \right)^{\fr{1}{2}} \right\} 
\le \fr{1}{8} \E \left\{\sup\limits_{0 \le t \le T} \|q^{\epsilon}\|_{L^2}^2 \right\} + C\|\tilde{g}\|_{L^2}^2 T
\end{align} 
This gives \eqref{SS1}.

\textbf{Proof of (ii).} Applying It\^o's lemma to the process $F(X_t (w))$ where $X_t(w) = \|q^{\epsilon}(t,w)\|_{L^2}^2$ obeys \eqref{ql2} and $F(\xi) = \xi^{\fr{p}{2}}$, we derive the energy equality 
\beg{align}
\d (\|q^{\epsilon}\|_{L^2}^2)^{\fr{p}{2}} 
&= - p \|q^{\epsilon}\|_{L^2}^{p-2} \|\l^{\fr{1}{2}} q^{\epsilon}\|_{L^2}^2 dt \nonumber
\\&\quad+ p \|q^{\epsilon}\|_{L^2}^{p-2} (\Delta \Phi, q^{\epsilon})_{L^2} dt
+ \fr{p}{2} \|q^{\epsilon}\|_{L^2}^{p-2} \sum\limits_{l=1}^{n} \|J_{\epsilon} \tilde{g}_l\|_{L^2}^2 dt \nonumber
\\&\quad +p \left(\fr{p}{2} -1 \right) \|q^{\epsilon}\|_{L^2}^{p-4} \sum\limits_{l=1}^{n} |(J_{\epsilon} \tilde{g}_l, q^{\epsilon})_{L^2} |^2 dt
+ \sum\limits_{l=1}^{n} p\|q^{\epsilon}\|_{L^2}^{p-2} (J_{\epsilon} \tilde{g}_l, q^{\epsilon})_{L^2} dW_l, 
\end{align} which yields the differential inequality
\beg{align}
&\d \|q^{\epsilon}\|_{L^2}^p + p \|q^{\epsilon}\|_{L^2}^{p-2} \|\l^{\fr{1}{2}} q^{\epsilon}\|_{L^2}^2 dt
\le p \|q^{\epsilon}\|_{L^2}^{p-1} \|\Delta \Phi\|_{L^2}dt \nonumber
\\&\quad\quad+ \fr{p}{2} (p-1) \|q^{\epsilon}\|_{L^2}^{p-2} \|\tilde{g}\|_{L^2}^2 dt
+  \sum\limits_{l=1}^{n} p\|q^{\epsilon}\|_{L^2}^{p-2} (J_{\epsilon} \tilde{g}_l, q^{\epsilon})_{L^2} dW_l. 
\end{align} 
In view of the bound
\be 
\|q^{\epsilon}\|_{L^2} \le \|\l^{\fr{1}{2}} q^{\epsilon}\|_{L^2},
\ee we have
\beg{align}  \la{qpeq}
&\d \|q^{\epsilon}\|_{L^2}^p  + \fr{p}{4} \|q^{\epsilon}\|_{L^2}^p dt + \fr{p}{2} \|q^{\epsilon}\|_{L^2}^{p-2} \|\l^{\fr{1}{2}} q^{\epsilon}\|_{L^2}^2 dt \nonumber
\\&\quad\le C_p \left(\|\Delta \Phi\|_{L^2}^p  + \|\tilde{g}\|_{L^2}^{p} \right)dt
+ \sum\limits_{l=1}^{n} p\|q^{\epsilon}\|_{L^2}^{p-2} (J_{\epsilon} \tilde{g}_l, q^{\epsilon})_{L^2} dW_l
\end{align} where we used Young's inequality to estimate 
\be 
p \|q^{\epsilon}\|_{L^2}^{p-1} \|\Delta \Phi\|_{L^2} 
\le C_p \|\Delta \Phi\|_{L^2}^p + \fr{p}{8} \|q^{\epsilon}\|_{L^2}^p
\ee and 
\be 
\fr{p}{2} (p-1) \|q^{\epsilon}\|_{L^2}^{p-2} \|\tilde{g}\|_{L^2}^2
\le C_p \|\tilde{g}\|_{L^2}^{p} + \fr{p}{8} \|q^{\epsilon}\|_{L^2}^p.
\ee
Integrating in time \eqref{qpeq} from $0$ to $t$ and taking the supremum over $[0,T]$, we obtain 
\beg{align}
&\sup\limits_{0 \le t \le T} \|q^{\epsilon}\|_{L^2}^p + \fr{p}{2} \int_{0}^{T} \|q^{\epsilon}\|_{L^2}^{p-2} \|\l^{\fr{1}{2}} q^{\epsilon}\|_{L^2}^2 ds \nonumber
\\&\quad\le 2\|q_0\|_{L^2}^p  
+ C_p \left(\|\Delta \Phi\|_{L^2}^p  + \|\tilde{g}\|_{L^2}^p \right) T 
+ 2\sup\limits_{0 \le t \le T} \left|\int_{0}^{t} \sum\limits_{l=1}^{n} p \|q^{\epsilon}\|_{L^2}^{p-2} (J_{\epsilon}\tilde{g}_l, q^{\epsilon})_{L^2} dW_l \right|.
\end{align} We estimate 
\beg{align}
&\E \left\{\sup\limits_{0 \le t \le T} \left|\int_{0}^{t} 2p\|q^{\epsilon}\|_{L^2}^{p-2} \sum\limits_{l=1}^{n} (J_{\epsilon} \tilde{g}_l, q^{\epsilon})_{L^2} dW_l \right| \right\}
\le C_p \E \left\{\left(\int_{0}^{T} \sum\limits_{l=1}^{n} \|q^{\epsilon}\|_{L^2}^{2p-4} (J_{\epsilon}\tilde{g}_l, q^{\epsilon})_{L^2}^2 dt \right)^{\fr{1}{2}}\right\} \nonumber
\\&\quad\le C_p \E \left\{\left(\int_{0}^{T} \|q^{\epsilon}\|_{L^2}^{2p-2} \|\tilde{g}\|_{L^2}^2 dt \right)^{\fr{1}{2}} \right\} 
\le \E \left\{\left(\sup\limits_{0 \le t \le T} \|q^{\epsilon}\|_{L^2}^{p-1} \right) \left(C_p \int_{0}^{T} \|\tilde{g}\|_{L^2}^2 dt \right)^{\fr{1}{2}} \right\} \nonumber
\\&\quad\le \left( 1 - \fr{1}{p} \right) \E \left\{\sup\limits_{0 \le t \le T} \|q^{\epsilon}\|_{L^2}^p \right\} + C_p \|\tilde{g}\|_{L^2}^p T^{\fr{p}{2}}
\end{align} and we obtain \eqref{SS2}.

\textbf{Proof of (iii).} We apply It\^o's lemma pointwise in $x$ to the processes $F(u_1^{\epsilon}(w))$ and $F(u_2^{\epsilon}(w))$ where $F(\xi) = \xi^2$, we add the resulting equations, and we integrate in the space variable over the torus $\TT^2$. We obtain the energy equality  
\beg{align}
\d \|u^{\epsilon}\|_{L^2}^2
&= -2(-\Delta u^{\epsilon}, u^{\epsilon})_{L^2} dt 
-2 (u^{\epsilon} \cdot \na u^{\epsilon}, u^{\epsilon})_{L^2} dt
-2 (q^{\epsilon} Rq^{\epsilon}, u^{\epsilon})_{L^2}dt
-2 (q^{\epsilon} \na \Phi, u^{\epsilon})_{L^2} dt \nonumber
\\&\quad\quad+2 (f, u^{\epsilon})_{L^2} dt
+ \sum\limits_{l=1}^{n} \|J_{\epsilon} g_l\|_{L^2}^2 dt
+ 2 \sum\limits_{l=1}^{n} (J_{\epsilon} g_l, u^{\epsilon})_{L^2} dW_l,
\end{align} which implies
\beg{align}
&\d \|u^{\epsilon}\|_{L^2}^2 + 2\|\na u^{\epsilon}\|_{L^2}^2 dt  \nonumber
\\&\quad= -2 (q^{\epsilon} Rq^{\epsilon} + q^{\epsilon} \na  \Phi -f, u^{\epsilon})_{L^2}dt
+ \sum\limits_{l=1}^{n} \|J_{\epsilon} g_l\|_{L^2}^2 dt + 2 \sum\limits_{l=1}^{n} (J_{\epsilon} g_l, u^{\epsilon})_{L^2} dW_l,
\end{align} where we used the cancellation
\be 
(u^{\epsilon} \cdot \na u^{\epsilon}, u^{\epsilon})_{L^2} = 0
\ee due to the divergence-free condition satisfied by $u^{\epsilon}$. 
By Ladyzhenskaya's interpolation inequality 
\be 
\|u^{\epsilon}\|_{L^4} \le C\|u^{\epsilon}\|_{L^2} + C\|u^{\epsilon}\|_{L^2}^{\fr{1}{2}} \|\na u^{\epsilon}\|_{L^2}^{\fr{1}{2}},
\ee and the boundedness of the Riesz transforms in $L^4$, we estimate
\beg{align}
|(q^{\epsilon} Rq^{\epsilon}, u^{\epsilon})_{L^2}| 
&\le \|q^{\epsilon}\|_{L^2} \|Rq^{\epsilon}\|_{L^4} \|u^{\epsilon}\|_{L^4} 
\le C\|q^{\epsilon}\|_{L^2} \|q^{\epsilon}\|_{L^4}\left(\|u^{\epsilon}\|_{L^2} + \|u^{\epsilon}\|_{L^2}^{\fr{1}{2}} \|\na u^{\epsilon}\|_{L^2}^{\fr{1}{2}} \right) \nonumber 
\\&\le C\|q^{\epsilon}\|_{L^2}^2 \|q^{\epsilon}\|_{L^4}^2 + \fr{1}{2} \|u^{\epsilon}\|_{L^2}^2 + \fr{1}{2}\|\na u^{\epsilon}\|_{L^2}^2.
\end{align} We also estimate 
\be 
|(q^{\epsilon} \na \Phi, u^{\epsilon})_{L^2}| \le \fr{1}{2} \|u^{\epsilon}\|_{L^2}^2 + \fr{1}{2} \|\na \Phi\|_{L^4}^2 \|q^{\epsilon}\|_{L^4}^2 
\ee and 
\be 
|(f, u^{\epsilon})_{L^2}| \le \fr{1}{2} \|u^{\epsilon}\|_{L^2}^2 + \fr{1}{2} \|f\|_{L^2}^2
\ee using H\"older's inequality followed by Young's inequality. 
We obtain the differential inequality 
\beg{align} 
\d \|u^{\epsilon}\|_{L^2}^2 + \|\na u^{\epsilon}\|_{L^2}^2 dt
&\le 3\|u^{\epsilon}\|_{L^2}^2 dt 
+ \|f\|_{L^2}^2 dt 
+ C\|q^{\epsilon}\|_{L^2}^2 \|q^{\epsilon}\|_{L^4}^2 dt \nonumber
\\&\quad\quad+ C \|\na \Phi\|_{L^4}^2\|q^{\epsilon}\|_{L^4}^2 dt
+ \|g\|_{L^2}^2 dt
+ 2\sum\limits_{l=1}^{n} (J_{\epsilon} g_l, u^{\epsilon})_{L^2} dW_l,
\end{align} hence
\beg{align}
&\d \left\{e^{-3t} \|u^{\epsilon}\|_{L^2}^2 \right\}(s) 
= -3e^{-3s} \|u^{\epsilon}\|_{L^2}^2 ds + e^{-3s} \d \|u^{\epsilon} (s)\|_{L^2}^2  \nonumber
\\&\quad\le -e^{-3s}\|\na u^{\epsilon}\|_{L^2}^2  ds + e^{-3s}  \left\{  \|f\|_{L^2}^2 ds + C\|q^{\epsilon}\|_{L^2}^2 \|q^{\epsilon}\|_{L^4}^2 ds + C \|\na \Phi\|_{L^4}^2\|q^{\epsilon}\|_{L^4}^2 ds \right\} \nonumber
\\&\quad\quad+ e^{-3s}\|g\|_{L^2}^2 ds + 2e^{-3s} \sum\limits_{l=1}^{n} (J_{\epsilon} g_l, u^{\epsilon})_{L^2} dW_l
\end{align} for all $s \in [0,t]$. Integrating in time from $0$ to $t$, we obtain 
\beg{align} \la{S3}
&e^{-3t} \|u^{\epsilon}(t)\|_{L^2}^2 + \int_{0}^{t} e^{-3s} \|\na u^{\epsilon}(s)\|_{L^2}^2 ds 
\le \|u_0\|_{L^2}^2  
+ \left(\|f\|_{L^2}^2 + \|g\|_{L^2}^2\right)t   \nonumber
\\&\quad\quad+ C \int_{0}^{t} \|q^{\epsilon}(s)\|_{L^2}^2 \|q^{\epsilon}(s)\|_{L^4}^2 ds 
+ C \int_{0}^{t} \|\na \Phi\|_{L^4}^2\|q^{\epsilon}(s)\|_{L^4}^2 ds \nonumber
\\&\quad\quad+ 2 \int_{0}^{t} e^{-3s} \sum\limits_{l=1}^{n} (J_{\epsilon} g_l, u^{\epsilon})_{L^2} dW_l(s).
\end{align}
We take the supremum in time over $[0,T]$ and apply $\E$. Using the continuous Sobolev embedding
\be 
H^{\fr{1}{2}} (\TT^2) \subset L^4(\TT^2)
\ee and \eqref{SS2} with $p=4$, we have 
\be 
\E \left\{\int_{0}^{T} \|q^{\epsilon}(s)\|_{L^2}^2 \|q^{\epsilon}(s)\|_{L^4}^2 ds  \right\} 
\le C\|q_0\|_{L^2}^4 + C \left(\|\Delta \Phi\|_{L^2}^4 + \|\tilde{g}\|_{L^2}^{2} \right)T + C\|\tilde{g}\|_{L^2}^2 T^2
\ee for all $t \in [0,T]$. From \eqref{SS1}, we have 
\be 
\E \left\{\int_{0}^{T} \|\na \Phi\|_{L^4}^2\|q^{\epsilon}(s)\|_{L^4}^2 ds \right\}
\le C\|\na \Phi\|_{L^4}^2 \left( \|q_0\|_{L^2}^2 + \|\l^{\fr{3}{2}} \Phi\|_{L^2}^2T + \|\tilde{g}\|_{L^2}^2 T \right)
\ee for all $t \in [0,T]$. We estimate 
\beg{align}
&\E \left\{\sup\limits_{0 \le t \le T}  \left|\int_{0}^{t} \sum\limits_{l=1}^{n} 2e^{-3s} (J_{\epsilon} g_l, u^{\epsilon})_{L^2} dW_l \right| \right\}
\le  \E \left\{\sup\limits_{0 \le t \le T} \left(e^{-\fr{3}{2}t} \|u^{\epsilon}(t)\|_{L^2}\right) \left(\int_{0}^{T} Ce^{-3t} \|g\|_{L^2}^2 dt \right)^{\fr{1}{2}} \right\} \nonumber
\\&\quad\le \fr{1}{2} \E \left\{ \sup\limits_{0 \le t \le T} \left(e^{-3t}  \|u^{\epsilon}(t)\|_{L^2}^2 \right) \right\} + C\|g\|_{L^2}^2 
\end{align} and we obtain \eqref{SS3}.

\textbf{Proof of (iv).} We apply It\^o's lemma pointwise in $x$ to the stochastic process $F(X_t(w))$ where $X_t(w) = q^{\epsilon}$ and $F(\xi) = \xi^4$. We have
\beg{align}
\d |q^{\epsilon}|^4 
&= -4 (q^{\epsilon})^3 u^{\epsilon} \cdot \na q^{\epsilon} dt
- 4 (q^{\epsilon})^3 \l q^{\epsilon} dt
+ 4 (q^{\epsilon})^3 \Delta \Phi dt \nonumber
\\&\quad\quad+ 6 \sum\limits_{l=1}^{n} (q^{\epsilon})^2 (J_{\epsilon}\tilde{g}_l)^2 dt
+ 4(q^{\epsilon})^3 \sum\limits_{l=1}^{n} J_{\epsilon} \tilde{g}_l d W_l.  
\end{align}
Integrating in the space over $\TT^2$, we obtain the energy equality
\beg{align}  \la{1}
\d \|q^{\epsilon}\|_{L^4}^4 
&= -4 (u^{\epsilon} \cdot \na q^{\epsilon}, (q^{\epsilon})^3)_{L^2} dt
- 4 (\l q^{\epsilon}, (q^{\epsilon})^3)_{L^2} dt
+ 4 (\Delta \Phi, (q^{\epsilon})^3)_{L^2} dt \nonumber
\\&+ 6(\sum\limits_{l=1}^{n} (J_{\epsilon} \tilde{g}_l)^2, (q^{\epsilon})^2)_{L^2} dt
+ 4 \sum\limits_{l=1}^{n} (J_{\epsilon} \tilde{g}_l, (q^{\epsilon})^3)_{L^2} dW_l.
\end{align} We note that 
\be \la{2}
(u^{\epsilon} \cdot \na q^{\epsilon}, (q^{\epsilon})^3)_{L^2} = 0
\ee due to the divergence-free condition for $u^{\epsilon}$. 
By the nonlinear Poincar\'e inequality for the fractional Laplacian in $L^4$ applied to the mean zero function $q^{\epsilon}$ (see \cite{AI,CGV}), we have 
\be  \la{3}
\int_{\TT^2} (q^{\epsilon})^3 \l q^{\epsilon} dx \ge c \|q^{\epsilon}\|_{L^4}^4.
\ee Using H\"older's inequality with exponents $4, 4/3$ and Young's inequality with exponents $4, 4/3$, we get 
\be  \la{4}
4|(\Delta \Phi, (q^{\epsilon})^3)_{L^2}|
\le 4\|\Delta \Phi\|_{L^{4}} \|(q^{\epsilon})^3\|_{L^{4/3}}
= 4\|\Delta \Phi\|_{L^4} \|q^{\epsilon}\|_{L^4}^3
\le c\|q^{\epsilon}\|_{L^4}^4 + C\|\Delta \Phi\|_{L^4}^4.
\ee We also bound
\be  \la{5}
6|(\sum\limits_{l=1}^{n} \left(J_{\epsilon} \tilde{g}_l)^2, (q^{\epsilon})^2 \right)_{L^2}|
\le 6\|q^{\epsilon}\|_{L^4}^2 \sum\limits_{l=1}^{n} \|\tilde{g}_l\|_{L^4}^2
\le c\|q^{\epsilon}\|_{L^4}^4 + C \left(\sum\limits_{l=1}^{n} \|\tilde{g}_l\|_{L^4}^2 \right)^{2},
\ee using H\"older and Young inequalities. 
Putting \eqref{1}--\eqref{5} together, we obtain the differential inequality 
\be \la{S4}
\d \|q^{\epsilon}\|_{L^4}^4 + c \|q^{\epsilon}\|_{L^4}^4 dt
\le C\|\Delta \Phi\|_{L^4}^4 dt + C \left(\sum\limits_{l=1}^{n} \|\tilde{g}_l\|_{L^4}^2 \right)^2 dt
+ 4 \sum\limits_{l=1}^{n} (J_{\epsilon} \tilde{g}_l, (q^{\epsilon})^3)_{L^2} dW_l. 
\ee Consequently,
\be 
\|q^{\epsilon}(t)\|_{L^4}^4 
+ c\int_{0}^{t} \|q^{\epsilon}\|_{L^4}^4 ds 
\le 2\|q_0\|_{L^4}^4  
+ C\|\Delta \Phi\|_{L^4}^4 t
+ C \left(\sum\limits_{l=1}^{n} \|\tilde{g}_l\|_{L^4}^2 \right)^2 t
+ 4 \int_{0}^{t} \sum\limits_{l=1}^{n} (J_{\epsilon} \tilde{g}_l, (q^{\epsilon})^3)_{L^2} dW_l
\ee for all $t \in [0,T]$. We take the supremum over $[0,T]$ and then we apply $\E$. We estimate 
\beg{align} \la{martingaleterm2}
&\E \left\{\sup\limits_{0 \le t \le T} \left|8\int_{0}^{t} \sum\limits_{l=1}^{n} (J_{\epsilon} \tilde{g}_l, (q^{\epsilon})^3)_{L^2} d W_l \right| \right\}
\le C\E \left\{\left(\int_{0}^{T} \sum\limits_{l=1}^{n} (J_{\epsilon}\tilde{g}_l, (q^{\epsilon})^3)_{L^2}^2 dt \right)^{\fr{1}{2}} \right\} \nonumber
\\&\quad\le C\E \left\{\left(\int_{0}^{T} \sum\limits_{l=1}^{n} \|\tilde{g}_l \|_{L^4}^2 \|(q^{\epsilon})^3\|_{L^{4/3}}^2 dt \right)^{\fr{1}{2}} \right\}
\le \E \left\{\sup\limits_{0 \le t \le T} \|q^{\epsilon}\|_{L^4}^3 \left(C\int_{0}^{T} \sum\limits_{l=1}^{n} \|\tilde{g}_l\|_{L^4}^2 dt \right)^{\fr{1}{2}} \right\} \nonumber
\\&\quad\le \fr{3}{4} \E \left\{\sup\limits_{0 \le t \le T} \|q^{\epsilon}\|_{L^4}^4 \right\} + C  \left(\sum\limits_{l=1}^{n} \|\tilde{g}_l \|_{L^4}^2 \right)^2  T^2
\end{align}
and we obtain \eqref{SS4}.

\textbf{Proof of (v).} We apply It\^o's lemma to the process $F(X_t(w))$ where $X_t(w) = \|q^{\epsilon}(t,w)\|_{L^4}^4$ and the twice differentiable function $F(\xi) = \xi^{\fr{p}{4}}$. We obtain
\beg{align}  
\d \|q^{\epsilon}\|_{L^4}^p 
&= - p\|q^{\epsilon}\|_{L^4}^{p-4} (\l q^{\epsilon}, (q^{\epsilon})^3)_{L^2} dt
+ p\|q^{\epsilon}\|_{L^4}^{p-4}  (\Delta \Phi, (q^{\epsilon})^3)_{L^2} dt \nonumber
\\&\quad\quad+ \fr{3}{2} p \|q^{\epsilon}\|_{L^4}^{p-4}(\sum\limits_{l=1}^{n} (J_{\epsilon} \tilde{g}_l)^2, (q^{\epsilon})^2)_{L^2} dt 
+ 2p \left(\fr{p}{4} -1 \right)\|q^{\epsilon}\|_{L^4}^{p-8} \sum\limits_{l=1}^{n} (J_{\epsilon} \tilde{g}_l, (q^{\epsilon})^3)_{L^2}^2 dt \nonumber
\\&\quad\quad+ p \|q^{\epsilon}\|_{L^4}^{p-4} \sum\limits_{l=1}^{n} (J_{\epsilon} \tilde{g}_l, (q^{\epsilon})^3)_{L^2} dW_l.
\end{align}
By H\"older's inequality with exponents $4/3, 4$ and Young's inequality with exponents $p/(p-2), p/2$, we have 
\beg{align} 
&2p \left(\fr{p}{4} -1 \right)\|q^{\epsilon}\|_{L^4}^{p-8} \sum\limits_{l=1}^{n} (J_{\epsilon} \tilde{g}_l, (q^{\epsilon})^3)_{L^2}^2
\le 2p \left(\fr{p}{4} -1 \right)\|q^{\epsilon}\|_{L^4}^{p-8}  \|(q^{\epsilon})^3\|_{L^{4/3}}^2 \sum\limits_{l=1}^{n} \|\tilde{g}_l \|_{L^4}^2 \nonumber
\\&\quad= 2p \left(\fr{p}{4} -1 \right) \|q^{\epsilon}\|_{L^4}^{p-8} \|q^{\epsilon}\|_{L^4}^6 \sum\limits_{l=1}^{n} \|\tilde{g}_l\|_{L^4}^2 
\le \fr{cp}{8} \|q^{\epsilon}\|_{L^4}^{p} + C\left(\sum\limits_{l=1}^{n} \|\tilde{g}_l \|_{L^4}^2 \right)^{\fr{p}{2}}.
\end{align}
We obtain 
\be \la{S5}
\d \|q^{\epsilon}\|_{L^4}^p  + \fr{cp}{2} \|q^{\epsilon}\|_{L^4}^p  dt
\le C\|\Delta \Phi\|_{L^4}^p dt + C\left(\sum\limits_{l=1}^{n} \|\tilde{g}_l \|_{L^4}^2 \right)^{\fr{p}{2}}  dt
+ p \|q^{\epsilon}\|_{L^4}^{p-4} \sum\limits_{l=1}^{n} (J_{\epsilon} \tilde{g}_l, (q^{\epsilon})^3)_{L^2} dW_l.
\ee Integrating \eqref{S5} in time from $0$ to $t$, taking the supremum over $[0,T]$, applying $\E$, and estimating
\beg{align} 
&\E \left\{\sup\limits_{0 \le t \le T} 2p\left|\int_{0}^{t} \|q^{\epsilon}\|_{L^4}^{p-4} \sum\limits_{l=1}^{n} (J_{\epsilon} \tilde{g}_l, (q^{\epsilon})^3)_{L^2} dW_l \right| \right\}\nonumber\\ 
&\quad\le \left(1 - \fr{1}{p} \right) \E \left\{\sup\limits_{0 \le t \le T} \|q^{\epsilon}\|_{L^4}^p \right\} + C_p \left(\sum\limits_{l=1}^{n} \|\tilde{g}_l\|_{L^4}^2 \right)^{\fr{p}{2}} T^{\fr{p}{2}} 
\end{align} we obtain \eqref{SS5}.

\textbf{Proof of (vi).} Using It\^o's lemma, we derive the energy equality
\beg{align}
\d (\|u^{\epsilon}\|_{L^2}^2)^{\fr{p}{2}}
&= -p \|u^{\epsilon}\|_{L^2}^{p-2} \|\na u^{\epsilon}\|_{L^2}^2 dt
+ p \|u^{\epsilon}\|_{L^2}^{p-2} (-q^{\epsilon} Rq^{\epsilon} - q^{\epsilon} \na \Phi + f, u^{\epsilon})_{L^2}dt \nonumber
\\&\quad\quad+ \fr{p}{2}\|u^{\epsilon}\|_{L^2}^{p-2} \sum\limits_{l=1}^{n} \|J_{\epsilon}g_{l} \|_{L^2}^2 dt
+ p \left(\fr{p}{2} -1 \right) \|u^{\epsilon}\|_{L^2}^{p-4} \sum\limits_{l=1}^{n} |(J_{\epsilon} g_l, u^{\epsilon})_{L^2}|^2 dt \nonumber
\\&\quad\quad+ p\|u^{\epsilon}\|_{L^2}^{p-2} \sum\limits_{l=1}^{n} (J_{\epsilon} g_l, u^{\epsilon})_{L^2} dW_l.
\end{align}
By Young's inequality with exponents $p/(p-2)$ and $p/2$, 
\be 
\fr{p}{2} \|u^{\epsilon}\|_{L^2}^{p-2} \sum\limits_{l=1}^{n} \|J_{\epsilon} g_l \|_{L^2}^2
\le  \fr{1}{5} \|u^{\epsilon}\|_{L^2}^p + C_p \|g\|_{L^2}^p 
\ee and 
\beg{align} 
p\left(\fr{p}{2} -1 \right)\|u^{\epsilon}\|_{L^2}^{p-4} \sum\limits_{l=1}^{n} |(J_{\epsilon} g_l, u^{\epsilon})_{L^2}|^2
&\le p\left(\fr{p}{2} -1 \right)\|u^{\epsilon}\|_{L^2}^{p-4} \|u^{\epsilon}\|_{L^2}^2 \|g\|_{L^2}^2 \nonumber
\\&\le \fr{1}{5}  \|u^{\epsilon}\|_{L^2}^p + C_p \|g\|_{L^2}^p.
\end{align} Similarly, using Young's inequality with exponents $p/(p-1)$ and $p$, 
\be 
p\|u^{\epsilon}\|_{L^2}^{p-2} |(f, u^{\epsilon})_{L^2}| \le p\|u^{\epsilon}\|_{L^2}^{p-2} \|u^{\epsilon}\|_{L^2}\|f\|_{L^2}
\le C_p\|f\|_{L^2}^p + \fr{1}{5}  \|u^{\epsilon}\|_{L^2}^p
\ee and 
\beg{align} 
p\|u^{\epsilon}\|_{L^2}^{p-2}|(q^{\epsilon} \na  \Phi, u^{\epsilon})_{L^2}| 
&\le p\|u^{\epsilon}\|_{L^2}^{p-2} \|u^{\epsilon}\|_{L^2} \|q^{\epsilon}\|_{L^2} \|\na \Phi\|_{L^{\infty}} \nonumber
\\&\le C_p \|\na \Phi\|_{L^{\infty}}^p \|q^{\epsilon}\|_{L^2}^p + \fr{1}{5} \|u^{\epsilon}\|_{L^2}^p.
\end{align} 
By Ladyzhenskaya's interpolation inequality and the boundedness of the Riesz transforms in $L^4(\TT^2)$, we have 
\beg{align} 
&p\|u^{\epsilon}\|_{L^2}^{p-2} |(-q^{\epsilon} Rq^{\epsilon}, u^{\epsilon})_{L^2}|
\le C_p\|u^{\epsilon}\|_{L^2}^{p-2} \|u^{\epsilon}\|_{L^4}\|q^{\epsilon}\|_{L^2}\|q^{\epsilon}\|_{L^4} \nonumber
\\&\quad\le C_p\|u^{\epsilon}\|_{L^2}^{p-2} \left(\|u^{\epsilon}\|_{L^2} + \|u^{\epsilon}\|_{L^2}^{\fr{1}{2}}\|\na u^{\epsilon}\|_{L^2}^{\fr{1}{2}} \right) \|q^{\epsilon}\|_{L^2}\|q^{\epsilon}\|_{L^4} \nonumber
\\&\quad\le \|u^{\epsilon}\|_{L^2}^{p} + \fr{p}{2} \|u^{\epsilon}\|_{L^2}^{p-2}\|\na u^{\epsilon}\|_{L^2}^2 + C_p \|q^{\epsilon}\|_{L^2}^p \|q^{\epsilon}\|_{L^4}^p \nonumber
\\&\quad\le \fr{1}{5}  \|u^{\epsilon}\|_{L^2}^{p} + \fr{p}{2} \|u^{\epsilon}\|_{L^2}^{p-2}\|\na u^{\epsilon}\|_{L^2}^2 + C_p \|q^{\epsilon}\|_{L^2}^{2p} +  C_p \|q^{\epsilon}\|_{L^4}^{2p}.
\end{align} This yields the differential inequality
\beg{align} 
&\d \|u^{\epsilon}\|_{L^2}^{p} + \fr{p}{2} \|u^{\epsilon}\|_{L^2}^{p-2} \|\na u^{\epsilon}\|_{L^2}^2 dt
\le \|u^{\epsilon}\|_{L^2}^p dt 
+ C_p \|g\|_{L^2}^p dt 
+ C_p \|f\|_{L^2}^p dt \nonumber
\\&\quad\quad+ C_p \|\na \Phi\|_{L^{\infty}}^p \|q^{\epsilon}\|_{L^2}^p dt 
+ C_p \|q^{\epsilon}\|_{L^2}^{2p} dt + C_p \|q^{\epsilon}\|_{L^4}^{2p} dt
+ p\|u^{\epsilon}\|_{L^2}^{p-2} \sum\limits_{l=1}^{n} (J_{\epsilon} g_l, u^{\epsilon})_{L^2} dW_l
\end{align} and thus
\beg{align}
&\d \left\{e^{-t} \|u^{\epsilon}\|_{L^2}^p \right\} (s) + e^{-s} \|u^{\epsilon}\|_{L^2}^{p-2} \|\na u^{\epsilon}\|_{L^2}^2 ds \nonumber
\\&\quad\le e^{-s} \left\{C_p \|g\|_{L^2}^{p}ds
+ C_p \|f\|_{L^2}^p ds
+ C_p \|\na \Phi\|_{L^{\infty}}^p \|q^{\epsilon}\|_{L^2}^p ds
+ C_p \|q^{\epsilon}\|_{L^2}^{2p} ds + C_p \|q^{\epsilon}\|_{L^4}^{2p} ds \right\} \nonumber
\\&\quad\quad+ pe^{-s} \|u^{\epsilon}\|_{L^2}^{p-2} \sum\limits_{l=1}^{n} (J_{\epsilon} g_l, u^{\epsilon})_{L^2} dW_l.
\end{align}
We integrate in time from $0$ to $t$, take the supremum over $[0,T]$, and apply $E$. We obtain 
\beg{align} \la{mod1}
&\E \left\{\sup\limits_{0 \le t \le T} \left(e^{-t} \|u^{\epsilon}(t)\|_{L^2}^p \right) \right\} 
+ \E \left\{\int_{0}^{T} e^{-t} \|u^{\epsilon}\|_{L^2}^{p-2} \|\na u^{\epsilon}\|_{L^2}^2 dt  \right\}\nonumber
\\&\quad\le C_p \left(\|g\|_{L^2}^p + \|f\|_{L^2}^p \right)
+ C_p \|\na \Phi\|_{L^{\infty}}^p \E \left\{ \int_{0}^{T} \|q^{\epsilon}\|_{L^2}^p dt \right\}
+ C_p \E \left\{\int_{0}^{T} \|q^{\epsilon}\|_{L^2}^{2p} dt \right\} \nonumber
\\&\quad\quad+ C_p \E \left\{\int_{0}^{T} \|q^{\epsilon}\|_{L^4}^{2p} dt  \right\}
+ \sup\limits_{0 \le t \le T} \left|\int_{0}^{t} 2pe^{-s} \|u^{\epsilon}\|_{L^2}^{p-2} \sum\limits_{l=1}^{n} (J_{\epsilon} g_l, u^{\epsilon})_{L^2} dW_l \right|.
\end{align}
We estimate  
\beg{align} \la{mod2}
&\E \left\{\sup\limits_{0 \le t \le T}  \left|\int_{0}^{t} 2pe^{-s} \|u^{\epsilon}\|_{L^2}^{p-2} \sum\limits_{l=1}^{n} (J_{\epsilon}q_l, u^{\epsilon})_{L^2} dW_l(s) \right| \right\} \nonumber
\\&\quad\le \left(1 - \fr{1}{p} \right) \E \left\{\sup\limits_{0 \le t \le T} \left(e^{-t}\|u^{\epsilon}(t)\|_{L^2}^p  \right) \right\}
+ C_p \|g\|_{L^2}^p T^{\fr{p}{2}}.
\end{align}
Putting \eqref{mod1} and \eqref{mod2} together, and using \eqref{SS2} and \eqref{SS5}, we obtain \eqref{SS6}.

\textbf{Proof of (vii).} We write the equation satisfied by $\na u^{\epsilon}$, apply It\^o's lemma, and integrate in the space variable. We obtain the energy equality 
\beg{align}
&\d \|\na u^{\epsilon}\|_{L^2}^2  + 2\|\Delta u^{\epsilon}\|_{L^2}^2
= 2(u^{\epsilon} \cdot \na u^{\epsilon}, \Delta u^{\epsilon})_{L^2} dt
+2 (q^{\epsilon}Rq^{\epsilon}, \Delta u^{\epsilon})_{L^2} dt \nonumber
\\&\quad\quad+ 2 (q^{\epsilon} \na \Phi, \Delta u^{\epsilon})_{L^2} dt
-2 (f, \Delta u^{\epsilon})_{L^2} dt
+ \|J_{\epsilon} \na g \|_{L^2}^2 dt
- 2\sum\limits_{l} (J_{\epsilon}g, \Delta u^{\epsilon})_{L^2} dW_l.
\end{align} The nonlinear term for the velocity vanishes, that is  
\be 
(u^{\epsilon} \cdot \na u^{\epsilon}, \Delta u^{\epsilon})_{L^2} = 0,
\ee and using H\"older's inequality, we obtain 
\beg{align} \la{expna}
&\d \|\na u^{\epsilon}\|_{L^2}^2 + 2 \|\Delta u^{\epsilon}\|_{L^2}^2 dt
\le C\|q^{\epsilon}\|_{L^4}^2 \|\Delta u^{\epsilon}\|_{L^2} dt
+ 2\|\na \Phi\|_{L^{\infty}} \|q^{\epsilon}\|_{L^2} \|\Delta u^{\epsilon}\|_{L^2} dt \nonumber
\\&\quad\quad+ 2\|f\|_{L^2} \|\Delta u^{\epsilon}\|_{L^2} dt
+ \|\na g \|_{L^2}^2 dt
- 2 \sum\limits_{l} (J_{\epsilon}g, \Delta u^{\epsilon})_{L^2} dW_l.
\end{align}
An application of Young's inequality yields the differential inequality 
\beg{align} \la{naudiff}
&\d \|\na u^{\epsilon}\|_{L^2}^2 +  \|\Delta u^{\epsilon}\|_{L^2}^2 dt
\le C\|q^{\epsilon}\|_{L^4}^4  dt
+ C\|\na \Phi\|_{L^{\infty}}^2 \|q^{\epsilon}\|_{L^2}^2  dt \nonumber
\\&\quad\quad+ C\|f\|_{L^2}^2 dt
+ \|\na g\|_{L^2}^2 dt
- 2 \sum\limits_{l} (J_{\epsilon}g, \Delta u^{\epsilon})_{L^2} dW_l.
\end{align}
We integrate \eqref{naudiff} in time from $0$ to $t$, take the supremum in time, and then apply $\E$. We obtain
\beg{align} \la{naudiffmod}
&\E\left\{\sup\limits_{0 \le t \le T} \|\na u^{\epsilon}\|_{L^2}^2 \right\} 
+ \E \left\{\int_{0}^{T} \|\Delta u^{\epsilon}\|_{L^2}^2 dt \right\}
\le 2\|\na u_0\|_{L^2}^2
+ C\left(\|\na g\|_{L^2}^2 + \|f\|_{L^2}^2 \right)T \nonumber
\\&\quad+ C\E\left\{\int_{0}^{T} \|q^{\epsilon}\|_{L^4}^4 dt \right\}
+ C \|\na \Phi\|_{L^{\infty}}^2 \E \left\{\int_{0}^{T}  \|q^{\epsilon}\|_{L^2}^2  dt \right\}
+ \sup\limits_{0 \le t \le T} \left|\int_{0}^{t} 4\sum\limits_{l} (J_{\epsilon}g, \Delta u^{\epsilon})_{L^2} dW_l \right|.
\end{align} We estimate the martingale term 
\beg{align} \la{martingaleterm}
&\E \left\{\sup\limits_{0 \le t \le T} \left|4 \int_{0}^{t} \sum\limits_{l} (J_{\epsilon}g_l, \Delta u^{\epsilon})_{L^2} dW_l \right| \right\} 
\le \E \left\{4\left(\int_{0}^{T} \sum\limits_{l} (J_{\epsilon}g_l, \Delta u^{\epsilon})_{L^2}^2 dt \right)^{\fr{1}{2}} \right\} \nonumber
\\&\quad\le \E \left\{4\left(\int_{0}^{T} \|\na g\|_{L^2}^2 \|\na u^{\epsilon}\|_{L^2}^2 dt \right)^{\fr{1}{2}} \right\} 
\le \E \left\{4\sup\limits_{0 \le t \le T} \|\na u^{\epsilon}\|_{L^2} \left(\int_{0}^{T} \|\na g\|_{L^2}^2 dt \right)^{\fr{1}{2}} \right\} \nonumber
\\&\quad\le \fr{1}{2} \sup\limits_{0 \le t \le T} \E \left\{\|\na u^{\epsilon}\|_{L^2}^2 \right\} + C\|\na g\|_{L^2}^2 T.
\end{align}
Putting \eqref{naudiffmod} and \eqref{martingaleterm} together, and using \eqref{SS1} and \eqref{SS4}, we get \eqref{SS7}.

Now we prove the existence of solutions for the stochastic electroconvection model \eqref{stochastic}. The proof uses ideas from \cite{MS}, where the authors investigated and determined the limiting drift for the stochastic 2D Navier-Stokes equations.

\beg{thm} \la{thm2} Let $T > 0$. Let $q_0 \in L^4$ have mean zero over $\TT^2$, and let $u_0 \in H^1$ be divergence-free. Suppose $\tilde{g}_l \in L^4$ and $g_l \in H^1$ for all $l \in \left\{1, ..., n\right\}. $
Then there exists a pair $(q, u)$ such that 
\be 
u \in L^{2} (\Omega; L^{\infty}(0,T ; H^1(\TT^2))) \cap L^2(\Omega, L^2(0,T;H^2(\TT^2))),
\ee
\be 
q \in L^{4} (\Omega; L^{\infty} (0,T; L^4(\TT^2))) \cap L^2(\Omega; L^2(0,T; H^{1/2}(\TT^2))),
\ee 
\be 
\d(q, \xi)_{L^2}  + (u \cdot \na q,\xi)_{L^2} dt + (\l q,\xi)_{L^2} dt = (\Delta \Phi,\xi)_{L^2} dt + \sum\limits_{l=1}^{n} (\tilde{g}_l, \xi)_{L^2} dW_l
\ee for any $\xi \in H^1(\TT^2)$ and a.e. $w \in \Omega$, and \be 
\d(u,v)_{L^2} + (u \cdot \na u + qRq, v)_{L^2} dt - (\Delta u,v)_{L^2} dt = (- q \na \Phi,v)_{L^2} dt + (f,v)_{L^2} dt + \sum\limits_{l=1}^{n} (g_l,v)_{L^2} dW_l
\ee  for any $v \in H^1(\TT^2)$ and a.e. $w \in \Omega$.
\end{thm}

\textbf{Proof:}  Let 
\be 
\mathcal{F}_1 (q^{\epsilon}, u^{\epsilon}) = u^{\epsilon} \cdot \na q^{\epsilon}
\ee and 
\be 
\mathcal{F}_2 (q^{\epsilon}, u^{\epsilon}) = u^{\epsilon} \cdot \na u^{\epsilon} + q^{\epsilon}Rq^{\epsilon}.
\ee We note that 
\beg{align} 
\|\mathcal{F}_1 \|_{H^{-1}}^2 
&\le \|u^{\epsilon}\|_{L^4}^2 \|q^{\epsilon}\|_{L^4}^2
\le C\left(\|u^{\epsilon}\|_{L^2}^2 + \|u^{\epsilon}\|_{L^2} \|\na u^{\epsilon}\|_{L^2}\right) \|q^{\epsilon}\|_{L^4}^2 \nonumber
\\&\le C \|u^{\epsilon}\|_{L^2}^4 + C\| q^{\epsilon}\|_{L^4}^4 + C \|u^{\epsilon}\|_{L^2}^2 \|\na u^{\epsilon}\|_{L^2}^2
\end{align} using Ladyzhenskaya's interpolation inequality, and 
\beg{align}
\|\mathcal{F}_2\|_{H^{-1}}^2 
&\le \|u^{\epsilon}\|_{L^4}^4 + \|q^{\epsilon}\|_{L^4}^2\|R q^{\epsilon}\|_{L^2}^2
\le  C\|u^{\epsilon}\|_{L^2}^4 + C\|u^{\epsilon}\|_{L^2}^2 \|\na u^{\epsilon}\|_{L^2}^2 + C\|\l^{\fr{1}{2}} q^{\epsilon}\|_{L^2}^2\|q^{\epsilon}\|_{L^2}^2
\end{align} using the boundedness of the Riesz transforms in $L^2(\TT^2)$. 
As a consequence of Proposition~\ref{prop1}, $\mathcal{F}_1$ and $\mathcal{F}_2$ are uniformly bounded in $L^2(\Omega, L^2(0,T; H^{-1}(\TT^2)))$. Therefore, up to subsequences, $u^{\epsilon}$ converges weakly to some function $u$ in \be 
L^{2} (\Omega; L^{\infty}(0,T ; H^1(\TT^2))) \cap L^2(\Omega, L^2(0,T; H^2(\TT^2))),
\ee $q^{\epsilon}$ converges weakly to some function $q$ in 
\be L^{4} (\Omega; L^{\infty} (0,T; L^4(\TT^2))) \cap L^2(\Omega; L^2(0,T; H^{1/2}(\TT^2))),
\ee and $\mathcal{F}_1 (q^{\epsilon}, u^{\epsilon})$ and  $\mathcal{F}_2 (q^{\epsilon}, u^{\epsilon})$ converge weakly to some functions $F_1$ and $F_2$, respectively, in 
\be L^2(\Omega, L^2(0,T; H^{-1}(\TT^2))).
\ee
Now we write the equations satisfied by $(q^{\epsilon}, u^{\epsilon})$ and $(q,u)$ as
\be  \la{thm21}
\d (q^{\epsilon}, u^{\epsilon}) + \mathcal{F}(q^{\epsilon}, u^{\epsilon}) dt  + (0,  \na p^{\epsilon} ) dt= (J_{\epsilon}\tilde{g}, J_{\epsilon}g) dW
\ee where 
\be 
\mathcal{F}(q^{\epsilon}, u^{\epsilon}) = (u^{\epsilon} \cdot \na q^{\epsilon} + \l q^{\epsilon} - \Delta \Phi, u^{\epsilon} \cdot \na u^{\epsilon} - \Delta u^{\epsilon} + q^{\epsilon}Rq^{\epsilon} + q^{\epsilon}\na \Phi - f),
\ee and 
\be  \la{thm22}
\d (q, u) + \mathcal{F}_0 dt = (\tilde{g}, g ) dW
\ee in $L^2(\Omega; L^2(0,T; H^{-1}(\TT^2)))$ where
\be 
\mathcal{F}_0 = (F_1 + \l q - \Delta \Phi, F_2 - \Delta u + q \nabla \Phi - f).
\ee We show that 
\be \la{dr3}
\mathcal{F}(q,u) = \mathcal{F}_0
\ee for almost every $w \in \Omega$.

We note that $(\l^{-1}q,u)$ obeys the energy equality 
\beg{align} \la{1111}
&\d \left(\|\l^{-\fr{1}{2}} q\|_{L^2}^2 + \|u\|_{L^2}^2 \right)
+ 2 (\mathcal{F}_0, (\l^{-1}q,u))_{L^2} dt \nonumber
\\&\quad= (\|\l^{-\fr{1}{2}}\tilde{g} \|_{L^2}^2 + \|g\|_{L^2}^2) dt
+ 2((\tilde{g}, g), (\l^{-1}q, u))_{L^2} dW.
\end{align}
We take a pair 
\be 
(\tilde{q}, \tilde{u}) \in L^4(\Omega; L^4(0,T; L^4(\TT^2))) \oplus L^2(\Omega; L^2(0,T; H^2)),
\ee where $\tilde{q}$ has mean zero and $\tilde{u}$ is divergence-free, and we define 
\be  \la{rcond}
r(t,w) = C_0 \int_{0}^{t} \left[\|\na \Phi\|_{L^{\infty}}^2 + \|\na \tilde{u}\|_{L^2}^2 + \|\na \tilde{u}\|_{L^2} + \|\tilde{q}\|_{L^4}^2 + \|\tilde{q}\|_{L^4}^4  +  \|\Delta \tilde{u}\|_{L^2}^2  \right] ds
\ee where $C_0$ is a large enough constant,  to be determined later.

The drift identification claim \eqref{dr3} is equivalent to showing that
\beg{align} \la{dr4}
\E \left\{\int_{0}^{T} 2e^{-r(t)} (\mathcal{F}(q,u) -  \mathcal{F}_0, (\l^{-1}\Psi_1, \Psi_2))_{L^2} dt \right\} \ge 0
\end{align} for all $(\Psi_1, \Psi_2) \in L^4(\Omega; L^4(0,T; L^4(\TT^2))) \oplus L^2(\Omega; L^2(0,T; H^2))$ such that $\Psi_1$ has mean zero and $\Psi_2$ is divergence-free.
Accordingly, we proceed to prove \eqref{dr4}.

Denoting $\d r(t)$ by $\dot{r}(t)$, we have 
\beg{align}
&\E \left\{\d \left[e^{-r(t)} \left(\|\l^{-\fr{1}{2}} q\|_{L^2}^2 + \|u\|_{L^2}^2 \right)  \right] + e^{-r(t)} (2\mathcal{F}_0 + \dot{r} (q, u)  , (\l^{-1}q,u) )_{L^2} dt \right\} \nonumber
\\&\quad= \E \left\{e^{-r(t)} \left(\|\l^{-\fr{1}{2}}\tilde{g}\|_{L^2}^2 +  \|g\|_{L^2}^2 \right) \right\}
\end{align} in view of \eqref{1111}, and consequently
\beg{align}
&\E \left\{-\int_{0}^{T} e^{-r(t)}(2\mathcal{F}_0 + \dot{r} (q,u), (\l^{-1}q, u))_{L^2} dt \right\} \nonumber
\\&\quad= \E \left\{e^{-r(T)}  \left(\|\l^{-\fr{1}{2}} q(T)\|_{L^2}^2 + \|u(T)\|_{L^2}^2 \right) -  \left(\|\l^{-\fr{1}{2}} q_0\|_{L^2}^2 + \|u_0\|_{L^2}^2 \right)\right\} \nonumber
\\&\quad\quad+ \E\left\{ - \int_{0}^{T} e^{-r(t)}  \left(\|\l^{-\fr{1}{2}}\tilde{g}\|_{L^2}^2 + \|g\|_{L^2}^2 \right) dt  \right\} \nonumber
\\&\quad\le \liminf_{\epsilon \to 0} \E \left\{e^{-r(T)}  \left(\|\l^{-\fr{1}{2}} q^{\epsilon}(T)\|_{L^2}^2 + \|u^{\epsilon}(T)\|_{L^2}^2 \right)\right\} + \lim \E \left\{-  \left(\|\l^{-\fr{1}{2}} J_{\epsilon}q_0\|_{L^2}^2 + \|J_{\epsilon}u_0\|_{L^2}^2 \right)  \right\} \nonumber
\\&\quad\quad+ \lim_{\epsilon \to 0} \E \left\{ - \int_{0}^{T} e^{-r(t)}  \left(\|\l^{-\fr{1}{2}}J_{\epsilon}\tilde{g}\|_{L^2}^2 + \|J_{\epsilon}g\|_{L^2}^2\right) dt   \right\} \nonumber
\\&\quad= \liminf\limits_{\epsilon \to 0} \E \left\{-\int_{0}^{T} e^{-r(t)}(2\mathcal{F}(q^{\epsilon}, u^{\epsilon}) + \dot{r} (q^{\epsilon},u^{\epsilon}), (\l^{-1}q^{\epsilon}, u^{\epsilon}))_{L^2} dt \right\}, 
\end{align} which implies that
\beg{align} \la{thm23}
&\E \left\{\int_{0}^{T} e^{-r(t)}(2\mathcal{F}_0 + \dot{r} (q,u), (\l^{-1}q, u))_{L^2} dt \right\} \nonumber
\\&\quad\ge \limsup_{\epsilon \to 0} \E \left\{\int_{0}^{T} e^{-r(t)}(2\mathcal{F}(q^{\epsilon}, u^{\epsilon}) + \dot{r} (q^{\epsilon},u^{\epsilon}), (\l^{-1}q^{\epsilon}, u^{\epsilon}))_{L^2} dt \right\}. 
\end{align}
We claim that 
\beg{align} \la{thm24}
&\E \left\{\int_{0}^{T} e^{-r(t)} (2\mathcal{F}(\tilde{q}, \tilde{u}) + \dot{r} (\tilde{q}, \tilde{u}), (\l^{-1}\tilde{q}, \tilde{u}) - (\l^{-1}q^{\epsilon}, u^{\epsilon}))_{L^2} dt \right\} \nonumber
\\&\quad\ge \E \left\{\int_{0}^{T} e^{-r(t)} (2\mathcal{F}(q^{\epsilon}, u^{\epsilon}) + \dot{r} (q^{\epsilon}, u^{\epsilon}), (\l^{-1} \tilde{q}, \tilde{u}) - (\l^{-1}q^{\epsilon}, u^{\epsilon}) )_{L^2} dt \right\}
\end{align}  for any $(\tilde{q}, \tilde{u}) \in L^4(\Omega; L^4(0,T; L^4(\TT^2))) \oplus L^2(\Omega; L^2(0,T; H^2))$ such that $\tilde{q}$ has mean zero and $\tilde{u}$ is divergence-free.

Suppose for now that the claim is true. Putting \eqref{thm23} and \eqref{thm24} together, we obtain 
\beg{align}
&\E \left\{\int_{0}^{T} e^{-r(t)} (2\mathcal{F} (\tilde{q}, \tilde{u}) + \dot{r} (\tilde{q}, \tilde{u}), (\l^{-1}\tilde{q}, \tilde{u})- (\l^{-1}q, u))_{L^2} dt \right\} \nonumber
\\&\quad= \lim_{\epsilon \to 0} \E\left\{\int_{0}^{T} e^{-r(t)} (2\mathcal{F} (\tilde{q}, \tilde{u}) + \dot{r} (\tilde{q}, \tilde{u}), (\l^{-1}\tilde{q}, \tilde{u})- (\l^{-1}q^{\epsilon}, u^{\epsilon}))_{L^2} dt \right\} \nonumber
\\&\quad\ge \liminf_{\epsilon \to 0} \E \left\{\int_{0}^{T} e^{-r(t)} (2\mathcal{F}(q^{\epsilon}, u^{\epsilon}) + \dot{r} (q^{\epsilon}, u^{\epsilon}), (\l^{-1} \tilde{q}, \tilde{u}) - (\l^{-1}q^{\epsilon}, u^{\epsilon}) )_{L^2} dt \right\} \nonumber
\\&\quad=  \E \left\{\int_{0}^{T} e^{-r(t)} (2\mathcal{F}_0 + \dot{r} (q, u), (\l^{-1} \tilde{q}, \tilde{u}) )_{L^2} dt \right\} \nonumber
\\&\quad\quad- \limsup_{\epsilon \to 0} \E \left\{\int_{0}^{T} e^{-r(t)} (2\mathcal{F}(q^{\epsilon}, u^{\epsilon}) + \dot{r} (q^{\epsilon}, u^{\epsilon}),  (\l^{-1}q^{\epsilon}, u^{\epsilon}) )_{L^2} dt \right\} \nonumber
\\&\quad\ge \E \left\{\int_{0}^{T} e^{-r(t)} (2\mathcal{F}_0 + \dot{r} (q,u), (\l^{-1} \tilde{q}, \tilde{u}) - (\l^{-1}q, u))_{L^2} dt \right\}
\end{align} for any $(\tilde{q}, \tilde{u}) \in L^4(\Omega; L^4(0,T; L^4(\TT^2))) \oplus L^2(\Omega; L^2(0,T; H^2))$ such that $\tilde{q}$ has mean zero and $\tilde{u}$ is divergence-free. 
Letting
\be 
(\tilde{q}, \tilde{u}) = (q, u) + \lambda \Psi
\ee where $\lambda > 0$ and $\Psi = (\Psi_1, \Psi_2) \in L^4(\Omega; L^4(0,T; L^4(\TT^2))) \oplus L^2(\Omega; L^2(0,T; H^2))$, $\Psi_1$ having mean zero and $\Psi_2$ being divergence-free, we obtain 
\beg{align}
&\E \left\{\int_{0}^{T} e^{-r(t)} (2\mathcal{F} ((q,u) + \lambda \Psi) + \dot{r} ((q,u) + \lambda \Psi ), \lambda (\l^{-1} \Psi_1, \Psi_2))_{L^2} dt \right\} \nonumber
\\&\quad\ge \E \left\{\int_{0}^{T} e^{-r(t)} (2\mathcal{F}_0 + \dot{r} (q,u),\lambda (\l^{-1} \Psi_1, \Psi_2))_{L^2} dt \right\}.
\end{align} We divide by $\lambda$, and then take the limit as $\lambda$ goes to zero. We obtain \eqref{dr4} from which we conclude that $\mathcal{F}_0 = \mathcal{F} (q,u)$. 

Now, in order to prove the claim \eqref{thm24}, it is enough to show that
\be \la{contidea}
(\mathcal{F}(\tilde{q}, \tilde{u}) - \mathcal{F} (q^{\epsilon}, u^{\epsilon}), (\l^{-1}(\tilde{q} - q^{\epsilon}) , \tilde{u} - u^{\epsilon}))_{L^2} + \dot{r} \left(\|\l^{-\fr{1}{2}} (\tilde{q} - q^{\epsilon}) \|_{L^2}^2 + \|\tilde{u} - u^{\epsilon} \|_{L^2}^2 \right) \ge 0.
\ee
Indeed, 
\beg{align}
&(\mathcal{F}(\tilde{q}, \tilde{u}) - \mathcal{F} (q^{\epsilon}, u^{\epsilon}), (\l^{-1}(\tilde{q} - q^{\epsilon}) , \tilde{u} - u^{\epsilon}))_{L^2} \nonumber
\\&= \int_{\TT^2} (\tilde{u} \cdot \na \tilde{q} - u^{\epsilon} \cdot \na q^{\epsilon}) \l^{-1} (\tilde{q} - q^{\epsilon}) 
+ \int_{\TT^2} \l (\tilde{q} - q^{\epsilon}) \l^{-1} (\tilde{q} - q^{\epsilon})
+ \int_{\TT^2} (\tilde{u} \cdot \na \tilde{u} - u^{\epsilon} \cdot \na u^{\epsilon}) \cdot (\tilde{u} - u^{\epsilon}) \nonumber
\\&\quad\quad- \int_{\TT^2} \Delta (\tilde{u} - u^{\epsilon}) \cdot (\tilde{u} - u^{\epsilon})
+ \int_{\TT^2} (\tilde{q} R\tilde{q} - q^{\epsilon} Rq^{\epsilon}) \cdot (\tilde{u} - u^{\epsilon})
+ \int_{\TT^2} (\tilde{q} - q^{\epsilon}) \na \Phi \cdot (\tilde{u} - u^{\epsilon}).
\end{align}
Integrating by parts, we have 
\be 
\int_{\TT^2} \l (\tilde{q} - q^{\epsilon}) \l^{-1} (\tilde{q} - q^{\epsilon}) 
-  \int_{\TT^2} \Delta (\tilde{u} - u^{\epsilon}) \cdot (\tilde{u} - u^{\epsilon})
= \|\tilde{q} - q^{\epsilon} \|_{L^2}^2 + \|\na (\tilde{u} - u^{\epsilon})\|_{L^2}^2.
\ee
By H\"older and Young inequalities, we have
\be 
\left|\int_{\TT^2} (\tilde{q} - q^{\epsilon}) \na \Phi \cdot (\tilde{u} - u^{\epsilon}) \right|
\le C\|\na \Phi\|_{L^{\infty}}^2\|\tilde{u} - u^{\epsilon}\|_{L^2}^2 + \fr{1}{4} \|\tilde{q} - q^{\epsilon}\|_{L^2}^2.
\ee We note that
\beg{align}
\int_{\TT^2} (\tilde{u} \cdot \na \tilde{u} - u^{\epsilon} \cdot \na u^{\epsilon}) \cdot (\tilde{u} - u^{\epsilon})  \nonumber
&= \int_{\TT^2} ((\tilde{u} - u^{\epsilon}) \cdot \na \tilde{u} ) \cdot (\tilde{u} - u^{\epsilon}) + \int_{\TT^2} (u^{\epsilon} \cdot \na (\tilde{u} - u^{\epsilon})) \cdot (\tilde{u} - u^{\epsilon})  
\\&=  \int_{\TT^2} ((\tilde{u} - u^{\epsilon}) \cdot \na \tilde{u} ) \cdot (\tilde{u} - u^{\epsilon})
\end{align} in view of the divergence-free condition satisfied by $u^{\epsilon}$, and hence 
\beg{align}
&\left| \int_{\TT^2} (\tilde{u} \cdot \na \tilde{u} - u^{\epsilon} \cdot \na u^{\epsilon}) \cdot (\tilde{u} - u^{\epsilon}) \right|
\le \|\na \tilde{u}\|_{L^2} \|\tilde{u} - u^{\epsilon}\|_{L^4}^2 \nonumber
\\&\quad\le C\|\na \tilde{u}\|_{L^2} \|\tilde{u} - u^{\epsilon}\|_{L^2} \|\na (\tilde{u} - u^{\epsilon})\|_{L^2}  + C\|\na \tilde{u}\|_{L^2} \|\tilde{u} - u^{\epsilon}\|_{L^2}^2 \nonumber
\\&\quad\le C\left(\|\na \tilde{u}\|_{L^2}^2 + \|\na \tilde{u}\|_{L^2} \right) \|\tilde{u} - u^{\epsilon}\|_{L^2}^2 
+ \fr{1}{4} \|\na (\tilde{u} - u^{\epsilon})\|_{L^2}^2
\end{align} where we used Ladyzhenskaya's interpolation inequality applied to $\tilde{u} - u^{\epsilon}$. 
Now, we write
\beg{align}\la{EQ1}
&\int_{\TT^2} (\tilde{u} \cdot \na \tilde{q} - u^{\epsilon} \cdot \na q^{\epsilon}) \l^{-1} (\tilde{q} - q^{\epsilon}) 
= \int_{\TT^2} ((\tilde{u} - u^{\epsilon}) \cdot \na \tilde{q}) \l^{-1} (\tilde{q} - q^{\epsilon}) \nonumber
\\&\quad\quad+ \int_{\TT^2} ((u^{\epsilon} - \tilde{u}) \cdot \na (\tilde{q} - q^{\epsilon})) \l^{-1} (\tilde{q} - q^{\epsilon})
+ \int_{\TT^2} (\tilde{u} \cdot \na (\tilde{q} - q^{\epsilon})) \l^{-1}(\tilde{q} - q^{\epsilon})
\end{align} and 
\beg{align}\la{EQ2}
&\int_{\TT^2} (\tilde{q} R\tilde{q} - q^{\epsilon} Rq^{\epsilon}) \cdot (\tilde{u} - u^{\epsilon})
= \int_{\TT^2} (\tilde{q} - q^{\epsilon}) R \tilde{q} \cdot (\tilde{u} - u^{\epsilon}) \nonumber
\\&\quad\quad+ \int_{\TT^2} (q^{\epsilon} - \tilde{q})R(\tilde{q} - q^{\epsilon}) \cdot (\tilde{u} - u^{\epsilon})
+ \int_{\TT^2} \tilde{q} R (\tilde{q} - q^{\epsilon}) \cdot (\tilde{u} - u^{\epsilon}). 
\end{align} Adding \eqref{EQ1} and \eqref{EQ2}, four terms cancel out, namely
\be 
\int_{\TT^2} ((u^{\epsilon} - \tilde{u}) \cdot \na (\tilde{q} - q^{\epsilon})) \l^{-1} (\tilde{q} - q^{\epsilon})
= - \int_{\TT^2} (q^{\epsilon} - \tilde{q})R(\tilde{q} - q^{\epsilon}) \cdot (\tilde{u} - u^{\epsilon})
\ee and 
\be 
\int_{\TT^2} ((\tilde{u} - u^{\epsilon}) \cdot \na \tilde{q}) \l^{-1} (\tilde{q} - q^{\epsilon})
= - \int_{\TT^2} \tilde{q} R (\tilde{q} - q^{\epsilon}) \cdot (\tilde{u} - u^{\epsilon}), 
\ee
due to the divergence-free condition satisfied by $u^{\epsilon} - \tilde{u}$. 
We estimate 
\beg{align}
&\left| \int_{\TT^2} (\tilde{q} - q^{\epsilon}) R \tilde{q} \cdot (\tilde{u} - u^{\epsilon}) \right|
\le \|R\tilde{q}\|_{L^4} \|\tilde{q} - q^{\epsilon}\|_{L^2} \|\tilde{u} - u^{\epsilon}\|_{L^4} \nonumber
\\&\quad\le C \|\tilde{q}\|_{L^4} \|\tilde{q} - q^{\epsilon}\|_{L^2} \left(\|\tilde{u} - u^{\epsilon}\|_{L^2} + \|\tilde{u} - u^{\epsilon}\|_{L^2}^{\fr{1}{2}} \|\na (\tilde{u} - u^{\epsilon}) \|_{L^2}^{\fr{1}{2}} \right) \nonumber
\\&\quad\le C\left(\|\tilde{q}\|_{L^4}^2 + \|\tilde{q}\|_{L^4}^4 \right) \|\tilde{u} - u^{\epsilon}\|_{L^2}^2 + \fr{1}{4} \|\tilde{q} - q^{\epsilon}\|_{L^2}^2 + \fr{1}{4} \|\na (\tilde{u} - u^{\epsilon}) \|_{L^2}^{2}
\end{align} using H\"older's inequality, the boundedness of the Riesz transforms in $L^4$, Ladyzhenskaya's inequality, and Young's inequality. In view of the commutator estimate (see \cite[Proposition 3]{AI}) 
\be 
\|[\l^{-\fr{1}{2}}, v \cdot \na ]\rho\|_{L^2} \le C \|\Delta v\|_{L^2} \|\rho\|_{L^2}
\ee
that holds for any divergence-free $v \in H^2$ and mean-zero $\rho \in L^2$,  we have 
\beg{align}
&\left|\int_{\TT^2} \tilde{u} \cdot \na (\tilde{q} - q^{\epsilon}) \l^{-1}(\tilde{q} - q^{\epsilon})\right|
= \left|\int_{\TT^2} \left[\l^{-\fr{1}{2}} (\tilde{u} \cdot \na (\tilde{q} - q^{\epsilon})) - \tilde{u} \cdot \na \l^{-\fr{1}{2}}(\tilde{q} - q^{\epsilon})  \right] \l^{-\fr{1}{2}}(\tilde{q} - q^{\epsilon})\right| \nonumber
\\&\quad\le C\|\Delta \tilde{u}\|_{L^2} \|\l^{-\fr{1}{2}} (\tilde{q} - q^{\epsilon})\|_{L^2} \|\tilde{q} - q^{\epsilon}\|_{L^2} 
\le C\|\Delta \tilde{u}\|_{L^2}^2 \|\l^{-\fr{1}{2}} (\tilde{q} - q^{\epsilon})\|_{L^2}^2 + \fr{1}{4} \|\tilde{q} - q^{\epsilon}\|_{L^2}^2.
\end{align} Therefore,
\beg{align}\la{EQ3}
&(\mathcal{F}(\tilde{q}, \tilde{u}) - \mathcal{F} (q^{\epsilon}, u^{\epsilon}), (\l^{-1}(\tilde{q} - q^{\epsilon}) , \tilde{u} - u^{\epsilon}))_{L^2}  \nonumber
\\&\quad+ C\left(\|\na \Phi\|_{L^{\infty}}^2 + \|\na \tilde{u}\|_{L^2}^2 + \|\na \tilde{u}\|_{L^2} + \|\tilde{q}\|_{L^4}^2 + \|\tilde{q}\|_{L^4}^4  +  \|\Delta \tilde{u}\|_{L^2}^2 \right)\left(\|\tilde{u} - u^{\epsilon}\|_{L^2}^2 + \|\l^{-\fr{1}{2}} (\tilde{q} - q^{\epsilon}) \|_{L^2}^2 \right) \nonumber
\\&\ge \fr{1}{4}\left(\|\na (\tilde{u} - u^{\epsilon})\|_{L^2}^2 + \|\tilde{q} - q^{\epsilon}\|_{L^2}^2 \right)
\ge 0.
\end{align} We choose the constant $C_0$ in \eqref{rcond} such that $C_0\ge C$, where $C$ is the absolute constant on the second line of inequality \eqref{EQ3}. Therefore, we obtain \eqref{thm24} and the proof of Theorem \ref{thm2} is complete.

\beg{rem}
Uniqueness of solutions is obtained as for the deterministic system \cite[Theorem 2]{AI}. Indeed, if we suppose the existence of two different solutions, and we write the equations obeyed by their difference, then we obtain deterministic equations which are independent of the noise. 
\end{rem}

\section{Electroconvection Semigroup} \la{sec3}

For each $t \ge 0$, we define
\be 
\mathcal{F}_t = \sigma \left(W_s: s \le t \right),
\ee that is, $\mathcal{F}_t$ is the smallest $\sigma$-algebra for which $W_s$ is measurable for all $s \le t$. 
Let $\tau$ be the stopping time random variable with respect to $\mathcal{F}_t$. 
For $u_0 = u_{\tau}(x,w)$ and $q_0 = q_{\tau}(x,w)$, we consider the electroconvection model \eqref{stochastic} in its variational form 
\be
\begin{cases} \la{semigroup1}
(q(t), \xi)_{L^2} + \int_{\tau}^{t} (u \cdot \na q (s), \xi)_{L^2} ds 
+ \int_{\tau}^{t} (\l q(s), \xi)_{L^2} ds
\\\quad= (q_{\tau}, \xi)_{L^2} + \int_{\tau}^{t} (\Delta \Phi, \xi)_{L^2} ds 
+ \sum\limits_{l=1}^{n} \int_{\tau}^{t} (\tilde{g}_l, \xi) dW_l (s)
\\(u(t), v)_{L^2} + \int_{\tau}^{t} (u \cdot \na u (s), v)_{L^2} ds 
\\\quad= (u_{\tau}, v)_{L^2} 
+ \int_{\tau}^{t} (-qRq (s) - q \nabla \Phi (s) + f, v)_{L^2} ds
+ \sum\limits_{l=1}^{n} \int_{\tau}^{t} (g_l, v)_{L^2} dW_{l}(s)
\end{cases}
\ee for any stopping time $\tau \le t \le T$, $\xi \in H^1(\TT^2)$ and $v \in H^1(\TT^2)$.  

\begin{thm} \la{semi2}
Let $\tau$ be a stopping time with respect to $\mathcal{F}_t$, and let $(q_{\tau}, u_{\tau})$ be $\mathcal{F}_{\tau}$ measurable random variables such that $0 \le \tau \le T$, $u_{\tau} \in L^2(\Omega; H^1(\TT^2))$ and $q_{\tau} \in L^4(\Omega; L^4(\TT^2))$. Suppose $\tilde{g}_l \in L^4$ and $g_l \in H^1$ for all $l \in \left\{1, ..., n \right\}$.
Then there exists a solution $(q,u)$ of \eqref{semigroup1} satisfying 
\be \la{unifbound1}
\E \left\{ \sup\limits_{\tau \le t \le T} \|q(t)\|_{L^4}^4  \right\} 
\le c_1 \E \left\{\|q_{\tau}\|_{L^4}^4 \right\} 
+ c_2 (\Phi, \tilde{g})T
+ c_3 (\tilde{g})T^2
\ee and
\be \la{unifbound2}
\E  \left\{\sup\limits_{\tau \le t \le T} \|\na u(t)\|_{L^2}^2 + \int_{\tau}^{T} \|\Delta u (s)\|_{L^2}^{2}  ds \right\} \le c_4 \E \left\{ \|\na u_{\tau}\|_{L^2}^2 + \|q_{\tau}\|_{L^4}^4   \right\} + c_5 (f, \Phi, g, \tilde{g})T + c_6(\tilde{g})T^2
\ee where $c_1$ and $c_4$ are universal constants, $c_2$ is a constant depending only on $\Phi$ and $\tilde{g}$, $c_3$  and $c_6$ are constants depending only on $\tilde{g}$, and $c_5$ is a constant depending on $f, \Phi, g$ and  $\tilde{g}$.
\end{thm}

\textbf{Proof:} The proof of \eqref{unifbound1} is similar to the proof of \eqref{SS4}. Indeed, we integrate \eqref{S4} from $\tau$ to $t$, we take the supremum over the time interval $[\tau, T]$, and then we apply $\E$. We estimate the noise term as in \eqref{martingaleterm2} and we obtain \eqref{unifbound1}.
As for the bound \eqref{unifbound2}, the proof is similar to the proof of \eqref{SS7}. Indeed, we integrate the differential inequality \eqref{naudiff} from $\tau$ to $t$, we take the supremum over $[\tau, T]$, and we take the expectation in $w$. We use \eqref{unifbound1} to estimate the charge density terms, and we bound the noise term as in \eqref{martingaleterm}. This gives \eqref{unifbound2}.  

\begin{thm} \la{Continuity} (Continuity) Let $(q_{\tau}^1, u_{\tau}^1)$ and $(q_{\tau}^2, u_{\tau}^2)$ be two initial data satisfying the assumptions of Theorem \ref{semi2}. Then the corresponding solutions $(q_1, u_1)$ and $(q_2, u_2)$ obey 
\begin{align}
&\|u_1(t) - u_2(t)\|_{L^2}^2 + \|\l^{-\fr{1}{2}} q_1(t) - \l^{-\fr{1}{2}}q_2(t)\|_{L^2}^2\nonumber\\
&\quad\le \exp \left\{CC(\tau, t)\right\} \left[\|u_{\tau}^1 - u_{\tau}^2\|_{L^2}^2 + \|\l^{-\fr{1}{2}} q_{\tau}^1 - \l^{-\fr{1}{2}} q_{\tau}^2 \|_{L^2}^2 \right]
\end{align} with probability 1, where 
\be 
C(\tau, t) = \int_{\tau}^{T} \left[\|\na \Phi\|_{L^{\infty}}^2 + \|\na u_1\|_{L^2}^2 + \|\na u_1\|_{L^2} + \|q_1\|_{L^4}^2 + \|q_1\|_{L^4}^4  +  \|\Delta u_1\|_{L^2}^2  \right] dt
\ee is well-defined and finite almost surely. 
\end{thm}

The proof is based on the same ideas used to prove \eqref{contidea}. We omit further details. 

\beg{prop} Let $(q_{\tau}, u_{\tau})$ be an initial data satisfying the conditions of Theorem \ref{semi2}.  Suppose $\tilde{g}_l \in L^4$ and $g_l \in H^1$ for all $l \in \left\{1, ..., n \right\}$. Then the unique solution $(q, u)$ of \eqref{semigroup1}  obeys
\be \la{weakbound}
\E \left\{ \sup\limits_{\tau \le t \le T} (\|\l^{-\fr{1}{2}} q\|_{L^2}^2 + \|u\|_{L^2}^2) \right\}
\le  \E \left\{\|\l^{-\fr{1}{2}} q_{\tau}\|_{L^2}^2 + \|u_{\tau}\|_{L^2}^2 + c_7(\Phi,f, g, \tilde{g}) \right\}e^{c_8(\Phi)T} 
\ee where $c_7$ is a positive constant depending only on $\Phi$, $f$, $g$, and $\tilde{g}$, and $c_8$ is a positive constant depending only on $\Phi$. 
\end{prop}

\textbf{Proof:} By It\^o's lemma, we have 
\beg{align} \la{EQ31}
\d \|\l^{-\fr{1}{2}}q\|_{L^2}^2 + 2 \|q\|_{L^2}^2 dt
&= -2 (u \cdot \na q, \l^{-1} q)_{L^2} dt + 2 (\Delta \Phi, \l^{-1}q)_{L^2} dt\nonumber
\\\quad\quad&+ \|\l^{-\fr{1}{2}} \tilde{g} \|_{L^2}^2  dt
+ 2 \sum\limits_{l=1}^{n} (\l^{-\fr{1}{2}}\tilde{g}_l, \l^{-\fr{1}{2}}q )_{L^2} dW_l
\end{align} and 
\beg{align}\la{EQ32}
\d \|u\|_{L^2}^2 + 2\|\na u \|_{L^2}^2 dt
&= - 2(u \cdot \na u, u)_{L^2} - 2(qRq, u)_{L^2} dt- 2(q\na \Phi, u)_{L^2}dt + 2(f, u)_{L^2}dt \nonumber
\\\quad\quad&+ \|g\|_{L^2}^2 dt + 2 \sum\limits_{l=1}^{n} (g_l, u)_{L^2} dW_l.
\end{align}
We add the equations \eqref{EQ31} and \eqref{EQ32}. Integrating by parts, we have
\be 
(u \cdot \na q, \l^{-1} q)_{L^2} = - (u \cdot Rq, q)_{L^2} = - (qRq, u)_{L^2},
\ee and using the cancellation 
\be 
(u \cdot \na u, u)_{L^2} = 0, 
\ee we obtain the differential equation 
\beg{align}\la{EQ33}
&\d \left\{\|\l^{-\fr{1}{2}}q\|_{L^2}^2 + \|u\|_{L^2}^2 \right\} + 2(\|q\|_{L^2}^2 + \|\na u\|_{L^2}^2) dt
= 2 (\Delta \Phi, \l^{-1}q)_{L^2} dt - 2(q\na \Phi, u)_{L^2} dt + 2(f, u)_{L^2} dt \nonumber
\\&\quad\quad+  \|\l^{-\fr{1}{2}} \tilde{g} \|_{L^2}^2 dt + \|g\|_{L^2}^2 dt
+ 2 \sum\limits_{l=1}^{n} (\l^{-\fr{1}{2}}\tilde{g}_l, \l^{-\fr{1}{2}}q )_{L^2} dW_l 
+ 2 \sum\limits_{l=1}^{n} (g_l, u)_{L^2} dW_l.
\end{align} From \eqref{EQ33}, we arrive at the differential inequality 
\beg{align}
&\d \left\{\|\l^{-\fr{1}{2}}q\|_{L^2}^2 + \|u\|_{L^2}^2 \right\} + (\|q\|_{L^2}^2 + \|\na u\|_{L^2}^2) dt
\le C(\|\l \Phi\|_{L^2}^2 + \|f\|_{L^2}^2) dt + C'(\|\na \Phi\|_{L^{\infty}}^2 + 1)\|u\|_{L^2}^2 dt \nonumber
\\&\quad\quad+  \|\l^{-\fr{1}{2}}\tilde{g}\|_{L^2}^2 dt + \|g\|_{L^2}^2 dt
+ 2 \sum\limits_{l=1}^{n} (\l^{-\fr{1}{2}}\tilde{g}_l, \l^{-\fr{1}{2}}q )_{L^2} dW_l 
+ 2 \sum\limits_{l=1}^{n} (g_l, u)_{L^2} dW_l.
\end{align}
Letting 
\be 
\rho = \|\na \Phi\|_{L^{\infty}}^2 + 1,
\ee we obtain 
\beg{align}
\d \left\{e^{-C'\rho t}(\|\l^{-\fr{1}{2}}q\|_{L^2}^2 + \|u\|_{L^2}^2) \right\} 
&\le C(\|\l \Phi\|_{L^2}^2 + \|f\|_{L^2}^2)e^{-C'\rho t} dt   
+  \|\l^{-\fr{1}{2}} \tilde{g} \|_{L^2}^2 dt   + \|g\|_{L^2}^2 dt \nonumber
\\&\quad\quad+ 2 \sum\limits_{l=1}^{n} (\l^{-\fr{1}{2}}\tilde{g}_l, \l^{-\fr{1}{2}}q )_{L^2} dW_l 
+ 2 \sum\limits_{l=1}^{n} (g_l, u)_{L^2} dW_l.
\end{align}
Integrating in time from $\tau$ to $t$, taking the supremum over $[\tau,T]$, applying the expectation $\E$ in $w$, and using martingale estimates, we obtain \eqref{weakbound}. 

We consider the space
\be 
\mathcal{H} = H^{-\fr{1}{2}} (\TT^2) \oplus L^2(\TT^2)
\ee consisting of vectors $(\xi, v)$ where $\xi \in H^{-\fr{1}{2}}$ has mean zero and $v \in L^2$ is divergence-free, and we consider the space
\be 
\mathcal{V} = L^4 (\TT^2)\oplus H^1(\TT^2) 
\ee
consisting of vectors $(\xi, v)$ where $\xi \in L^4$ has mean zero and $v \in H^1$ is divergence-free. We define the norms $\|\cdot\|_{\mathcal{H}}$ and $\|\cdot\|_{\mathcal{V}}$  by 
\be 
\|(\xi, v)\|_{\mathcal{H}}^2 = \|\l^{-\fr{1}{2}} \xi\|_{L^2}^2 + \|v\|_{L^2}^2.
\ee and 
\be 
\|(\xi, v)\|_{\mathcal{V}}^2 = \| \xi\|_{L^4}^2 + \|v\|_{H^1}^2
\ee respectively. 
Let $C_{g}^{0}(\mathcal{V}, \|\cdot\|_{\mathcal{H}})$ be the space of real continuous functions $h$ on the space $(\mathcal{V}, \|\cdot\|_{\mathcal{H}}) $, with growth 
\be  \la{grow}
|h(\xi, v)| \le C(1 + \|\l^{-\fr{1}{2}} \xi\|_{L^2}^2 + \|v\|_{L^2}^2).
\ee We point out that continuity of $h$ on the space $(\mathcal{V}, \|\cdot\|_{\mathcal{H}})$ means that if $(\xi_n ,v_n) \in \mathcal{V}$ converges to $(\xi, v)$ in the norm $\|\cdot\|_{\mathcal{H}}$, then $h(\xi_n, v_n)$ converges to $h(\xi, v)$.

Let $(\Phi(t,s), t \ge s \ge 0)$ be the semigroup associated to the electroconvection model \eqref{stochastic}
\be 
\Phi(t,s): C_{g}^{0} (\mathcal{V}) \rightarrow C_{g}^{0} (\mathcal{V})
\ee defined by 
\be 
\Phi (t,s)h(\xi, v) = \E \left\{h(q(t,s; \xi), u(t,s; v)) \right\}
\ee
where $(q(t,s; \xi), u(t,s; v))$ is the solution of \eqref{stochastic} with deterministic initial data $(q_s(x), u_s(x)) = (\xi(x), v(x))$. 

We note that the uniqueness of solutions in $\mathcal{V}$ (see \cite{AI}) imply that $(\Phi(t,s), t \ge s \ge 0)$  is indeed a semigroup. Moreover, $(\Phi(t,s), t \ge s \ge 0)$ is a $\mathcal{H}$-Markov Feller semigroup:

\beg{thm} \la{MFeller} ($\mathcal{H}$-Markov Feller Continuity) The semigroup $\Phi(t,s)$ is Markov-Feller on $C_g^0 (\mathcal{V},\|\cdot\|_{\mathcal{H}})$ in the sense that if $h \in C_g^{0}(\mathcal{V}, \|\cdot\|_{\mathcal{H}})$ and $\left\{(\xi_n, v_n) \right\}_{n=1}^{\infty}$ is a sequence in $\mathcal{V}$ converging to $(\xi,v) \in \mathcal{V}$ in the norm $\|\cdot\|_{\mathcal{H}}$, then 
\be 
\Phi(t,s) h(\xi_n, v_n) \rightarrow \Phi(t,s)h(\xi,v),
\ee and if $t_n \rightarrow s$, then 
\be 
\Phi(t_n, s) h(\xi,v) \rightarrow h(\xi,v)
\ee for any $(\xi, v) \in \mathcal{V}$. 
\end{thm}

\textbf{Proof:} Fix $h \in C_g^{0}(\mathcal{V}, \|\cdot\|_{\mathcal{H}})$. Suppose $(\xi_n, v_n)$ converges to $(\xi, v)$ in $(\mathcal{V}, \|\cdot\|_{\mathcal{H}})$, that is 
\be  \la{semi4}
\|\l^{-\fr{1}{2}} (\xi_n - \xi)\|_{L^2}^2 + \|v_n - v\|_{L^2} \rightarrow 0.
\ee In view of the continuity property given in Theorem~\ref{Continuity}, we have 
\be 
\|q(t,s; \xi_n) - q(t,s; \xi)\|_{H^{-\fr{1}{2}}} \rightarrow 0
\ee and 
\be 
\|u(t,s; v_n) - u(t,s;v)\|_{L^2} \rightarrow 0.
\ee Since $h$ is continuous on $(\mathcal{V}, \|\cdot\|_{\mathcal{H}})$, we conclude that 
\be 
h(q(t,s;\xi_n), u(t,s;v_n)) \rightarrow h(q(t,s; \xi), u(t,s;v))
\ee and hence 
\be 
\E \left\{h(q(t,s;\xi_n), u(t,s;v_n) \right\} \rightarrow \E \left\{h(q(t,s; \xi), u(t,s;v)) \right\}
\ee by the Lebesgue Dominated Convergence Theorem (which can be applied due to the growth condition \eqref{grow}, the bound \eqref{weakbound}, and the convergence \eqref{semi4}) yielding the boundedness of the sequence of initial datum $(\xi_n, v_n)$ in the $\mathcal{H}$-norm.   

Now, suppose that $\left\{t_n\right\}_{n=1}^{\infty}$ is a sequence of positive times converging to $s$, and $(\xi, v) \in \mathcal{V}$. Noting that the solution $(q(t,s;\xi), u(t,s;v))$ of \eqref{semigroup1} belongs to the space 
\be 
L^2(\Omega; C^0(s,T; H^{-\fr{1}{2}}(\TT^2))) \oplus L^2(\Omega; C^0(s,T; L^2(\TT^2))),
\ee we obtain 
\be
\E \left\{h(q(t_n, s ; \xi), u(t_n, s ;v))\right\} \rightarrow \E \left\{h(\xi,v)\right\}
\ee  due to the continuity of $h$ in $(\mathcal{V}, \|\cdot\|_{\mathcal{H}})$ and the Lebesgue Dominated Convergence Theorem. This ends the proof of Theorem \ref{MFeller}.

\section{Existence of an Invariant Measure in the Absence of Potential} \la{sec5}

In this section, we consider the electroconvection system 
\be \begin{cases} \la{invsys}
\d q + u \cdot \na q dt + \l q dt = \sum\limits_{l=1}^{n} \tilde{g}_l dW_l
\\\d u + u \cdot \na u dt - \Delta u dt + \na p dt = - q Rq dt + fdt + \sum\limits_{l=1}^{n} g dW_l
\\ \na \cdot u = 0
\end{cases}
\ee in  $\mathbb{T}^2 \times [0,T] \times \Omega$. We note that if the initial charge density and velocity are assumed to have a zero spatial average, then the solution $(q, u)$ will have mean zero over $\TT^2$ for all positive times $t \ge 0$.

Let $\dot{L}^p(\TT^2)$ and $\dot{H}^s(\TT^2)$ be the spaces of $L^p(\TT^2)$ and $H^s(\TT^2)$ functions with zero spatial averages respectively.
Let $H$ and $V$ be the spaces of $L^2(\TT^2)$ and $H^1(\TT^2)$ functions that are divergence-free and mean zero respectively.  
Let 
\be 
\dot{\mathcal{H}} = \dot{H}^{-\fr{1}{2}} (\TT^2) \oplus H
\ee and 
\be 
\dot{\mathcal{V}} = \dot{L}^4 (\TT^2) \oplus V
\ee  with 
\be 
\|(q, u)\|_{\dot{\mathcal{H}}}^2 = \|\l^{-\fr{1}{2}} q\|_{L^2}^2 + \|u\|_{L^2}^2
\ee and 
\be 
\|(q, u)\|_{\dot{\mathcal{V}} }^2 = \|q\|_{L^4}^2 + \|\na u\|_{L^2}^2
\ee respectively. 
See \cite{AI} for details on the notation and functional setting. We note that $\dot{\mathcal{V}}$ is compactly embedded in $\dot{\mathcal{H}}$.

We define the Markov transition kernels $\left\{P_t \right\}_{t \ge 0}$ associated to the electroconvection model \eqref{invsys} as 
\be 
P_t(q_0, u_0, A) = \mathbb{P} ((q(t, q_0), u(t, u_0)) \in A).
\ee
These kernels are defined on $\dot{\mathcal{V}}$ and are $\dot{\mathcal{H}}$-Feller as shown in Theorem \ref{MFeller}.

We will show that the solution $(q,u)$ of \eqref{invsys} lies in 
\be 
L^2(\Omega, L^2(0,T; {\dot{H}}^{\fr{3}{2}}(\TT^2))) \oplus L^2(\Omega, L^2(0,T; {H}^2(\TT^2) \cap H)
\ee and the bounds are linear in $T$, hence the Krylov-Bogoliubov procedure can be applied in order to prove the existence of an invariant measure.

The rigorous estimates in this section can be done by taking a viscous system approximating \eqref{invsys}, deriving the bounds for the mollified solution, and then inheriting them to the solution of \eqref{invsys} using the lower semi-continuity of the norms. We present formal proofs, omitting the
approximation. We need the following propositions:

\beg{prop} Let $p$ be an even integer such that $p \in \left\{4\right\} \cup [8, \infty)$. Let $q_0 \in \dot{L}^4$. Suppose $\tilde{g}_l \in \dot{L}^4$ for all $l \in \left\{1, ..., n\right\}$. Then there exist a positive constant $\Gamma_1$ depending only on $\|q_0\|_{L^4}$, $p$ and some universal constants, and a positive constant  $\Gamma_2$ depending only on $\tilde{g}$, $p$ and some universal constants such that 
\be \la{invariant1}
\int_{0}^{t} \E \|q(s)\|_{L^4}^p ds
\le \Gamma_1 (\|q_0\|_{L^4}) + \Gamma_2 (\tilde{g})t
\ee for all $t \ge 0$. Here $\Gamma_1 = 0$ if $q_0 = 0$.
\end{prop}

\textbf{Proof:} The $p$-th power of the $L^4$ norm of $q$ obeys the energy inequality
\be 
\d \|q\|_{L^4}^p  + \fr{cp}{2} \|q\|_{L^4}^p  dt
\le C\left(\sum\limits_{l=1}^{n} \|\tilde{g}_l \|_{L^4}^2 \right)^{\fr{p}{2}}  dt
+ p \|q\|_{L^4}^{p-4} \sum\limits_{l=1}^{n} (\tilde{g}_l, q^3)_{L^2} dW_l.
\ee Integrating in time from $0$ to $t$ and applying $\E$, we obtain the desired bound \eqref{invariant1}. 

\beg{prop} Let $u_0 \in V$ and $q_0 \in \dot{L}^4$. Suppose $g_l \in V$ and $\tilde{g}_l \in \dot{L}^4$ for all $l \in \left\{1, ..., n\right\}$. Then there exist positive constants $\Gamma_3, \Gamma_5$  depending only on $\|\na u_0\|_{L^2}$, $\|q_0\|_{L^4}$ and some universal constants, and  positive constants  $\Gamma_4, \Gamma_6$ depending only on $f, g$, $\tilde{g}$ and some universal constants such that 
\be \la{invariant2}
\E \|\na u(t)\|_{L^2}^2 + \E \left\{\int_{0}^{t} \|\Delta u(s)\|_{L^2}^2 ds \right\}
\le \Gamma_3 (\|\na u_0\|_{L^2}, \|q_0\|_{L^4}) + \Gamma_4 (f, g, \tilde{g})t,
\ee and 
\be \la{invariant3}
\E \left\{\int_{0}^{t} \|\na u(s)\|_{L^2}^2 \|\Delta u(s)\|_{L^2}^2 ds \right\}
\le \Gamma_5 (\|\na u_0\|_{L^2}, \|q_0\|_{L^4}) + \Gamma_6 (f, g, \tilde{g})t,
\ee  hold for all $t \ge 0$. Here $\Gamma_3 = \Gamma_5 = 0$ if $u_0 = q_0 = 0$.
\end{prop}

\textbf{Proof:} The $L^2$ norm of $\na u$ obeys
\beg{align} 
\d \|\na u \|_{L^2}^2  + 2\|\Delta u \|_{L^2}^2
= 2 (q Rq, \Delta u)_{L^2} dt 
-2 (f, \Delta u)_{L^2} dt
+ \|\na g \|_{L^2}^2 dt
- 2\sum\limits_{l} (g, \Delta u)_{L^2} dW_l.
\end{align} In view of H\"older's inequality, Young's inequality, and the boundedness of the Riesz transforms on $L^4(\TT^2)$, we get the energy inequality
\beg{align}
\d \|\na u \|_{L^2}^2  + \|\Delta u \|_{L^2}^2 dt
\le C\|q\|_{L^4}^4 dt + C\|f\|_{L^2}^2 dt + \|\na g\|_{L^2}^2 dt - 2\sum\limits_{l} (g, \Delta u)_{L^2} dW_l. 
\end{align} Integrating in time from $0$ to $t$ and applying $\E$, we obtain 
\beg{align} 
\E \|\na u(t)\|_{L^2}^2 
&+ \int_{0}^{t} \E \|\Delta u(s)\|_{L^2}^2 ds
\le \|\na u_0\|_{L^2}^2 \nonumber
\\&\quad + C \left(\|f\|_{L^2}^2 + \|\na g\|_{L^2}^2 \right)t
+ C \E \left\{\int_{0}^{t} \|q(s)\|_{L^4}^4 ds \right\}. 
\end{align}
In view of the bound \eqref{invariant1} applied with $p=4$, we obtain \eqref{invariant2}. 

By It\^o's lemma, we have
\beg{align} \la{evolprod}
&\d \|\na u\|_{L^2}^4
= -4 \|\na u \|_{L^2}^2 \|\Delta u\|_{L^2}^2 dt
+ 4\|\na u\|_{L^2}^2 (qRq  - f, \Delta u)_{L^2} dt \nonumber
\\&\quad\quad+2 \|\na u\|_{L^2}^2 \|\na g\|_{L^2} ^2dt
+ 4 \sum\limits_{l=1}^{n} |(g_l, \Delta u)_{L^2}|^2 dt
-4 \|\na u\|_{L^2}^2 \sum\limits_{l=1}^{n} (g_l, \Delta u)_{L^2} dW_l
\end{align} 
hence
\beg{align} \la{part4}
&\d \|\na u\|_{L^2}^4
+ 4\|\na u\|_{L^2}^2 \|\Delta u\|_{L^2}^2 dt
\le 4\|\na u\|_{L^2}^2 \|\Delta u\|_{L^2} \left(C\|q\|_{L^4}^2  + \|f\|_{L^2} \right)dt \nonumber
\\&\quad\quad+ 2\|\na u\|_{L^2}^2 \|\na g\|_{L^2}^2 dt
+ 4\|\na u\|_{L^2}^2 \sum\limits_{l=1}^{n} \|\na g_l\|_{L^2}^2 dt
 -4 \|\na u\|_{L^2}^2 \sum\limits_{l=1}^{n} (g_l, \Delta u)_{L^2} dW_l.
\end{align}
From \eqref{part4}, we obtain the differential inequality
\beg{align} 
&\d \|\na u\|_{L^2}^4 
+  \|\na u\|_{L^2}^4 dt
+  \|\na u\|_{L^2}^2\|\Delta u\|_{L^2}^2 dt
\le C\|q\|_{L^4}^4 \|\na u\|_{L^2}^2 dt \nonumber
\\&\quad\quad+ C(\|f\|_{L^2}^2 + \|\na g\|_{L^2}^2) \|\na u\|_{L^2}^2 dt
 -4 \|\na u\|_{L^2}^2 \sum\limits_{l=1}^{n} (g_l, \Delta u)_{L^2} dW_l,
\end{align} and by Young's inequality, we get
\beg{align} \la{EQ51}
&\d \|\na u\|_{L^2}^4 
+  \fr{1}{2} \|\na u\|_{L^2}^4 dt
+  \|\na u\|_{L^2}^2\|\Delta u\|_{L^2}^2 dt
\le C\|q\|_{L^4}^{8} dt \nonumber
\\&\quad\quad+ C(\|f\|_{L^2}^4 + \|\na g\|_{L^2}^4)  dt
 -4 \|\na u\|_{L^2}^2 \sum\limits_{l=1}^{n} (g_l, \Delta u)_{L^2} dW_l.
\end{align}
We integrate in time from $0$ to $t$ and we apply $\E$. In view of the bound \eqref{invariant1} applied with $p=8$, we obtain \eqref{invariant3}.

\beg{prop} Let $u_0 \in V$ and $q_0 \in \dot{L}^4$. Suppose $g_l \in V$ and $\tilde{g}_l \in \dot{L}^4$ for all $l \in \left\{1, ..., n\right\}$. Then there exist a positive constant $\Gamma_7$ depending only on $\|\na u_0\|_{L^2}$, $\|q_0\|_{L^4}, \tilde{g}, g, f$ and some universal constants, and a positive constant  $\Gamma_8$ depending only on $f, g$, $\tilde{g}$ and some universal constants such that 
\be \la{invariant4}
\E \left\{\int_{0}^{t} \|\na u(s)\|_{L^2}^2 \|\Delta u(s)\|_{L^2}^2 \|q(s)\|_{L^4}^4 ds \right\}
\le \Gamma_7 (\|\na u_0\|_{L^2}, \|q_0\|_{L^4}, \tilde{g}, g, f) + \Gamma_8 (f, g, \tilde{g})t
\ee holds for all $t \ge 0$. Here $\Gamma_7 = 0$ if $u_0 = q_0 = 0$.
\end{prop}

\textbf{Proof:} The stochastic process $\|\na u\|_{L^2}^4 \|q\|_{L^4}^4$ obeys
\be 
\d \left[\|\na u\|_{L^2}^4 \|q\|_{L^4}^4 \right]
= \|\na u\|_{L^2}^4 \d \|q\|_{L^4}^4 
+ \|q\|_{L^4}^4 \d \|\na u\|_{L^2}^4
+ \d \|\na u\|_{L^2}^4 \cdot \d \|q\|_{L^4}^4.
\ee The $4$-th power of the $L^2$ norm of $\na u$ evolves according to \eqref{evolprod} whereas the $4$-th power of the $L^4$ norm of $q$ evolves according to 
\beg{align}  
\d \|q\|_{L^4}^4 
= - 4 (\l q, q^3)_{L^2} dt
+ 6 \sum\limits_{l=1}^{n}  (\tilde{g}_l^2, q^2)_{L^2} dt
+ 4 \sum\limits_{l=1}^{n} (\tilde{g}_l, q^3)_{L^2} dW_l.
\end{align}
Consequently, the product $\|\na u\|_{L^2}^4 \|q\|_{L^4}^4$ satisfies the energy equality
\beg{align}
&\d \left[\|q\|_{L^4}^4 \|\na u\|_{L^2}^4 \right]
= -4 \|\na u\|_{L^2}^4 (\l q, q^3)_{L^2} dt
+ 6\|\na u\|_{L^2}^4 \sum\limits_{l=1}^{n} (\tilde{g}_l^2, q^2)_{L^2} dt \nonumber
\\&\quad + 4\|\na u\|_{L^2}^4 \sum\limits_{l=1}^{n} (\tilde{g}_l, q^3)_{L^2} dW_l 
-4\|q\|_{L^4}^4 \|\na u\|_{L^2}^2\|\Delta u\|_{L^2}^2 dt 
+ 4\|q\|_{L^4}^4 \|\na u\|_{L^2}^2 (qRq - f, \Delta u)_{L^2} dt \nonumber
\\&\quad\quad +2\|q\|_{L^4}^4 \|\na u\|_{L^2}^2 \|\na g\|_{L^2}^2 dt
+ 4\|q\|_{L^4}^4 \sum\limits_{l=1}^{n} (g_l, \Delta u)_{L^2}^2 dt
- 4\|q\|_{L^4}^4 \|\na u\|_{L^2}^2 \sum\limits_{l=1}^{n} (g_l, \Delta u)_{L^2} dW_l \nonumber
\\&\quad\quad\quad- 16\|\na u\|_{L^2}^2 \sum\limits_{l=1}^{n} (\tilde{g}_l, q^3)_{L^2} (g_l, \Delta u)_{L^2} dt
\end{align} which yields the energy inequality
\beg{align} \la{proof1}
&\d \left[\|q\|_{L^4}^4 \|\na u\|_{L^2}^4 \right]
+4c \|\na u\|_{L^2}^4 \|q\|_{L^4}^4 dt
+4\|q\|_{L^4}^4 \|\na u\|_{L^2}^2\|\Delta u\|_{L^2}^2 dt \nonumber
\\&\quad \le 6\|\na u\|_{L^2}^4 \sum\limits_{l=1}^{n} (\tilde{g}_l^2, q^2)_{L^2} dt
+ 4\|q\|_{L^4}^4 \|\na u\|_{L^2}^2 (qRq - f, \Delta u)_{L^2} dt \nonumber
\\&\quad\quad +2\|q\|_{L^4}^4 \|\na u\|_{L^2}^2 \|\na g\|_{L^2}^2 dt
+ 4\|q\|_{L^4}^4 \sum\limits_{l=1}^{n} (\na g_l, \na u)_{L^2}^2 dt
- 16\|\na u\|_{L^2}^2 \sum\limits_{l=1}^{n} (\tilde{g}_l, q^3)_{L^2} (g_l, \Delta u)_{L^2} dt \nonumber
\\&\quad\quad\quad - 4\|q\|_{L^4}^4 \|\na u\|_{L^2}^2 \sum\limits_{l=1}^{n} (g_l, \Delta u)_{L^2} dW_l 
+ 4\|\na u\|_{L^2}^4 \sum\limits_{l=1}^{n} (\tilde{g}_l, q^3)_{L^2} dW_l
\end{align} in view of the Poincar\'e inequality for the fractional Laplacian in $L^4$.
By the Cauchy-Schwartz inequality, Young's inequality and the Poincar\'e inequality applied to the mean zero function $\na u$, we estimate
\beg{align} 
&\left|6\|\na u\|_{L^2}^4 \sum\limits_{l=1}^{n} (\tilde{g}_l, q^2)_{L^2} \right|
\le 6\|\na u\|_{L^2}^4 \left(\sum\limits_{l=1}^{n} \|\tilde{g}_l\|_{L^4}^2 \right) \|q\|_{L^4}^2  \nonumber
\\&\quad \le \fr{c}{8} \|\na u\|_{L^2}^4 \|q\|_{L^4}^4 
+ C\left(\sum\limits_{l=1}^{n} \|\tilde{g}_l\|_{L^4}^2 \right)^2 \|\na u\|_{L^2}^2 \|\Delta u\|_{L^2}^2.
\end{align}
The boundedness of the Riesz transforms on $L^4(\TT^2)$ yields
\beg{align}
&\left|4\|q\|_{L^4}^4 \|\na u\|_{L^2}^2 (qRq - f, \Delta u)_{L^2} \right|
\le C\|q\|_{L^4}^6 \|\na u\|_{L^2}^2 \|\Delta u\|_{L^2} 
+ C\|q\|_{L^4}^4 \|\na u\|_{L^2}^2 \|\Delta u\|_{L^2} \|f\|_{L^2} \nonumber
\\&\quad \le \fr{1}{8} \|q\|_{L^4}^4\|\na u\|_{L^2}^2 \|\Delta u\|_{L^2}^2 
+ \fr{c}{8} \|q\|_{L^4}^4 \|\na u\|_{L^2}^4 
+ C\|q\|_{L^4}^{12}
+ C\|q\|_{L^4}^4 \|f\|_{L^2}^4. 
\end{align}
We bound 
\be 
2\|q\|_{L^4}^4 \|\na u\|_{L^2}^2 \|\na g\|_{L^2}^2 
\le \fr{c}{8} \|q\|_{L^4}^4 \|\na u\|_{L^2}^4 
+ C\|\na g\|_{L^2}^4 \|q\|_{L^4}^4
\ee and 
\be 
4\|q\|_{L^4}^4 \sum\limits_{l=1}^{n} (\na g_l, \na u)_{L^2}^2 
\le 4\|q\|_{L^4}^4\|\na u\|_{L^2}^2 \|\na g\|_{L^2}^2
\le \fr{c}{8} \|q\|_{L^4}^4 \|\na u\|_{L^2}^4 
+ C\|\na g\|_{L^2}^4 \|q\|_{L^4}^4
\ee using Young's inequality. 
Finally, we estimate
\beg{align} \la{proof2}
&\left|16\|\na u\|_{L^2}^2 \sum\limits_{l=1}^{n} (\tilde{g}_l, q^3)_{L^2} (g_l, \Delta u)_{L^2}\right|
\le 16\|\na u\|_{L^2}^3 \|q\|_{L^4}^3  \left(\sum\limits_{l=1}^{n} \|\tilde{g}_l\|_{L^4} \right) \left(\sum\limits_{l=1}^{n} \|\na g_l\|_{L^2} \right) \nonumber
\\&\quad \le \fr{c}{8} \|\na u\|_{L^2}^4 \|q\|_{L^4}^4 
+ C\left(\sum\limits_{l=1}^{n} \|\tilde{g}_l\|_{L^4} \right)^4 \left(\sum\limits_{l=1}^{n} \|\na g_l\|_{L^2} \right)^4.
\end{align}
Putting \eqref{proof1}--\eqref{proof2} together, we end up with the differential inequality
\beg{align}
&\d \left[\|q\|_{L^4}^4 \|\na u\|_{L^2}^4 \right]
+ c\|\na u\|_{L^2}^4 \|q\|_{L^4}^4 dt
+ \|q\|_{L^4}^4 \|\na u\|_{L^2}^2 \|\Delta u\|_{L^2}^2 dt \nonumber 
\\&\quad \le K_1(\tilde{g}) \|\na u\|_{L^2}^2 \|\Delta u\|_{L^2}^2 dt
+ K_2 (f, g) \|q\|_{L^4}^4 dt
+ K_3(g, \tilde{g})
+ C\|q\|_{L^4}^{12} dt
\end{align} where $K_1>0$ is a constant depending only on $\tilde{g}$, $K_2>0$ is a constant depending only on $f$ and $g$, $K_3>0$ is a constant depending only on $g$ and $\tilde{g}$, and $C$ is a positive universal constant.
We integrate in time from $0$ to $t$ and we apply $\E$. The bound \eqref{invariant1} applied with $p=4$ and $p=12$ together with the bound \eqref{invariant3} gives the desired estimate \eqref{invariant4}. 

\beg{prop} Let $u_0 \in V$ and $q_0 \in \dot{H}^1$. Suppose $g_l \in V$ and $\tilde{g}_l \in \dot{H}^1$ for all $l \in \left\{1, ..., n\right\}$. Then there exists a positive constant $\Gamma_9$ depending only on $\|\na u_0\|_{L^2}$, $\|\na q_0\|_{L^2}, \tilde{g}, g, f$ and some universal constants, and a positive constant  $\Gamma_{10}$ depending only on $f, g$, $\tilde{g}$ and some universal constants such that 
\be \la{invariant5}
\E \left\{\int_{0}^{t} \|\l^{\fr{3}{2}} q(s)\|_{L^2}^2 ds \right\} 
\le \Gamma_9 (\|\na u_0\|_{L^2}, \|\na q_0\|_{L^2}, \tilde{g}, g, f) + \Gamma_{10} (f, g, \tilde{g})t
\ee for any $t \ge 0$. Here $\Gamma_9 = 0$ if $u_0 = q_0 = 0$. 
\end{prop}

\textbf{Proof:} By It\^o's lemma, we have
\beg{align} 
&\d \|\na q\|_{L^2}^2
+ 2\|\l^{\fr{3}{2}} q\|_{L^2}^2 dt\nonumber\\
&\quad= 2(u \cdot \na q, \Delta q)_{L^2} dt
+ \|\na \tilde{g}\|_{L^2}^2 dt
- 2 \sum\limits_{l=1}^{n} (\tilde{g}_l, \Delta q)_{L^2} dW_l.
\end{align} 
We estimate the nonlinear term
\be 
|(u \cdot \na q, \Delta q)_{L^2}| \le \|\na u\|_{L^4} \|\na q\|_{L^{\fr{8}{3}}}^2 \le C\|\na u\|_{L^4} \|q\|_{L^4}^{\fr{1}{2}} \|\l^{\fr{3}{2}} q\|_{L^2}^{\fr{3}{2}} 
\ee 
using H\"older's inequality, and the interpolation inequality \cite[Proposition 2]{AI}
\be 
\|\l^{\fr{3}{2}} q\|_{L^2}^2 \ge C\|q\|_{L^4}^{-\fr{2}{3}} \|\na q\|_{L^{\fr{8}{3}}}^{\fr{8}{3}}.
\ee
We obtain the stochastic energy inequality
\be \la{EQ56}
\d \|\na q\|_{L^2}^2
+ \|\l^{\fr{3}{2}} q\|_{L^2}^2 dt
\le C\|\na u\|_{L^2}^{2} \|\Delta u\|_{L^2}^2 \|q\|_{L^4}^{2} dt
+ \|\na \tilde{g}\|_{L^2}^2 dt
- 2 \sum\limits_{l=1}^{n} (\tilde{g}_l, \Delta q)_{L^2} dW_l
\ee and by Young's inequality, we obtain
\beg{align} 
\d \|\na q\|_{L^2}^2
&+ \|\l^{\fr{3}{2}} q\|_{L^2}^2 dt
\le C\|\na u\|_{L^2}^{2} \|\Delta u\|_{L^2}^2 \|q\|_{L^4}^{4} dt \nonumber
\\&\quad+ C\|\na u\|_{L^2}^{2} \|\Delta u\|_{L^2}^2 dt
+ \|\na \tilde{g}\|_{L^2}^2 dt
- 2 \sum\limits_{l=1}^{n} (\tilde{g}_l, \Delta q)_{L^2} dW_l.
\end{align}
We integrate in time from $0$ to $t$ and we apply $\E$. In view of \eqref{invariant3} and \eqref{invariant4}, we obtain \eqref{invariant5}. 

The above propositions give

\beg{prop} \la{prop4} Suppose $g_l \in H^1(\TT^2)$ and $\tilde{g}_l \in H^1(\TT^2)$ for all $l \in \left\{1, ..., n\right\}$. Let 
\be 
\nu_T (A) = \fr{1}{T} \int_{0}^{T} \PP ((q(s), u(s)) \in A ) ds.
\ee Then $\left\{\nu_T \right\}$ is tight for $u_0 = q_0 = 0$. 
\end{prop}

\textbf{Proof:} Suppose $u_0 = q_0 = 0$. Using the bounds \eqref{invariant2} and \eqref{invariant5}, we have 
\be \la{prop41}
\E \int_{0}^{T} \|\Delta u\|_{L^2}^2 ds \le  \Gamma_4 (f, g, \tilde{g}) T
\ee and 
\be \la{prop42}
\E \int_{0}^{T} \|\l^{\fr{3}{2}} q\|_{L^2}^2 ds 
\le \Gamma_{10}(f, g, \tilde{g})T
\ee for all $T \ge 0$. 
Now, let $R > 0$, and let $B_R$ be the ball of radius $R$ in $\dot{H}^\fr{3}{2}(\TT^2) \oplus (\dot{H}^2(\TT^2) \cap H)$ which is compact in $\dot{\mathcal{V}}$. By Chebyshev's inequality, 
\beg{align} 
\sup\limits_{T > 0} \nu_{T} (B_R^c) 
&= \sup\limits_{T > 0} \fr{1}{T} \int_{0}^{T} \PP (\|(q, u)\|_{\dot{H}^{\fr{3}{2}}(\TT^2) \oplus (\dot{H}^2(\TT^2) \cap H)} \ge R) dt \nonumber
\\&\le \fr{1}{R^2} \sup\limits_{T > 0} \fr{1}{T} \int_{0}^{T} \E (\|(q, u)\|_{\dot{H}^\fr{3}{2}(\TT^2) \oplus (\dot{H}^2(\TT^2) \cap H)}^2) dt \rightarrow 0
\end{align} as $R \rightarrow \infty$ in view of the bounds \eqref{prop41} and \eqref{prop42} that are linear in $T$. Therefore, the family $\left\{\nu_T \right\}$ is tight, ending the proof of Proposition~\ref{prop4}.

As a consequence of the Krylov-Bogoliubov averaging procedure, we obtain 

\beg{thm} Suppose that $g_l \in H^1(\TT^2)$ and $\tilde{g}_l \in H^1(\TT^2)$ for all $l \in \left\{1, ..., n\right\}$. There exists an invariant measure for the Markov transition kernels associated with \eqref{invsys}.
\end{thm}

\section{Subcritical Case} \la{sec8}

For $\alpha > 1$, we consider the stochastic subcritical electroconvection model
\be \begin{cases} \la{stochasticsub}
\d q + u \cdot \na q dt + \l^{\alpha} q dt =  \sum\limits_{l=1}^{n} \tilde{g}_l dW_l 
\\\d u + u \cdot \na u dt - \Delta u dt + \na p dt = - q Rq dt  + fdt + \sum\limits_{l=1}^{n} g_l dW_l
\\ \na \cdot u = 0
\end{cases}
\ee on $\mathbb{T}^2 \times [0,T] \times \Omega$, with initial data $u(x,0) = u_0$ and $q(x,0) = q_0$.
Here $\l^{\alpha}$ is the fractional Laplacian of order $\alpha$. 

The existence and uniqueness of solutions is obtained as for the critical case (when the fractional Laplacian is of order 1). 

The solution $(q,u)$ of \eqref{stochasticsub} obeys
\be 
\int_{0}^{T} \E \|\l^{\fr{\alpha}{2}}q \|_{L^2}^2 dt \le \|q_0\|_{L^2}^2 + \Gamma_{11}(\tilde{g}) T
\ee and 
\be 
\int_{0}^{T} \E \|\Delta u\|_{L^2}^2 dt \le \Gamma_{12}(\|\na u_0\|_{L^2}, \|q_0\|_{L^4}) + \Gamma_{13}(f, g, \tilde{g})T
\ee for all $T \ge 0$, where $\Gamma_{11}$ is a positive constant depending only on $\tilde{g}$ and some universal constants, $\Gamma_{12}$ is a positive constant depending only on the initial data, and $\Gamma_{13}$ is a positive constant depending only on $f, g, \tilde{g}$ and some universal constants. 
In view of the compactness of $H^{\fr{\alpha}{2}}(\TT^2)$ in $L^4(\TT^2)$ for $\alpha > 1$ (which does not hold in the critical case), the Krylov-Bogoliubov averaging procedure implies automatically the existence of an invariant measure.    

\beg{thm} \la{subun} Suppose $g_l \in H^1(\TT^2)$ and $\tilde{g}_l \in L^4(\TT^2)$ for all $l \in \left\{1, ..., n\right\}$. Then there exists an invariant measure for the Markov transition kernels associated with \eqref{stochasticsub}. 
\end{thm}


{\bf {Acknowledgment.}} We thank N. Glatt-Holtz for suggesting to add stochastic forcing in electroconvection.

\end{document}